\documentclass[11pt,a4paper]{article}
\usepackage{ifthen}
\usepackage{amsmath}
\usepackage{amssymb}
\usepackage{latexsym}
\usepackage{bbm}
\usepackage{graphicx}
\usepackage{pdfsync}
\usepackage{verbatim}
\usepackage{fancyhdr}
\usepackage{url}               
\usepackage[english]{babel}
\usepackage[applemac]{inputenc}
\usepackage{hyperref}

\usepackage{amsthm}

\hypersetup{
	colorlinks=true,     
	citecolor=blue,
	linkcolor=red,     
	}

\newtheorem{theorem}{Theorem}[section]
\newtheorem{lemma}[theorem]{Lemma}
\newtheorem{assump}[theorem]{Assumptions}
\newtheorem{proposition}[theorem]{Proposition}

\theoremstyle{remark}
\newtheorem{remark}{Remark}

\numberwithin{equation}{section}
\numberwithin{figure}{section}

\newcommand{\bp}{{\it Proof. }}
\newcommand{\ep}{\hfill $\square$\\}

\newcommand{\bpl}{{\it Proof of Lemma }} 
\newcommand{\epl}{\hfill $\square$\\}

\newcommand{\bpp}{{\it Proof of Proposition }} 
\newcommand{\epp}{\hfill $\square$\\}

\newcommand{\be}{\begin{equation}}
\newcommand{\ee}{\end{equation}}
\newcommand{\bea}{\begin{eqnarray}}
\newcommand{\eea}{\end{eqnarray}}
\newcommand{\bee}{\begin{eqnarray*}}
\newcommand{\eee}{\end{eqnarray*}}

\def\pa{\partial}

\def\na{\nabla}

\def\RR{\mathbb{R}}

\def\ni{\noindent}

\def\eps{\vare}

\def\eps{\varepsilon}

\catcode`@=11
\def\supess{\mathop{\operator@font Sup\,ess}}
\catcode`@=12

\def\RR{\mathbb{R}}

\def\ds{\displaystyle}

\def\ni{\noindent}

\def\bar#1{{\overline #1}}

\def\pa{\partial}
\def\na{\nabla}

\newcommand{\ud}{\mathrm{d}}

\begin{document}

\renewcommand{\refname}{References}
\bibliographystyle{alpha}

\pagestyle{fancy}
\fancyhead[L]{ }
\fancyhead[R]{}
\fancyfoot[C]{}
\fancyfoot[L]{ }
\fancyfoot[R]{}
\renewcommand{\headrulewidth}{0pt} 
\renewcommand{\footrulewidth}{0pt}

\newcommand{\montitre}{On the fractional diffusion for the linear Boltzmann equation with drift and general cross-section}

\newcommand{\auteur}{\textsc{ Dahmane Dechicha
}}
\newcommand{\affiliation}{Laboratoire de recherche AGM, CY Cergy Paris Universit\'e. UMR CNRS 8088 \\
2 Avenue Adolphe Chauvin 95302 Cergy-Pontoise Cedex, France \\
Department of Mathematics, New York University in Abu Dhabi \\
Saadiyat Island, P.O. Box 129188, Abu Dhabi, UAE \\
\url{dahmane.dechicha@nyu.edu} }

 \begin{center}
{\bf  {\LARGE \montitre}}\\ \bigskip \bigskip
 {\large\auteur}\\ \bigskip \smallskip
 \affiliation \\ \bigskip
\today
 \end{center}
 \begin{abstract}
This paper is devoted to the hydrodynamic limit for the linear Boltzmann equation, in the case of a heavy tail equilibrium and a cross section which depends on the space variable and which degenerates for large velocities, without symmetry assumptions. For an appropriate time scale, a macroscopic equation with an elliptic operator, which is equivalent to the fractional Laplacian, is obtained. This problem has been addressed in \cite{MMM} for a space-independent cross section, using a Fourier-Laplace transform, where a fractional diffusion equation has been obtained, and revisited in \cite{M} for a space-dependent but a bounded cross section,  using the moments method by introducing an auxiliary problem. In this work, we will adapt the latter method to generalize both results.   \end{abstract}


\tableofcontents

 \pagestyle{fancy}
\fancyhead[R]{\thepage}
\fancyfoot[C]{}
\fancyfoot[L]{}
\fancyfoot[R]{}
\renewcommand{\headrulewidth}{0.2pt} 
\renewcommand{\footrulewidth}{0pt} 

\thispagestyle{empty}

\section{Introduction}
\subsection{Setting of the problem}
In this paper, we consider the following parameterized linear kinetic equation:
\begin{equation}\label{BL-eps}
\left\{\begin{array}{l}
\theta(\varepsilon) \partial_t f^\varepsilon + \varepsilon (v-j^\eps_F)\cdot \na_x f^\varepsilon = Q(f^\varepsilon), \hspace{0.5cm} x \in \RR^{d} , v \in \RR^{d} , t > 0 , \\
\\
    f^\varepsilon(0,x,v) = f_0(x,v), \hspace{3.2cm} x \in \RR^{d} ,   v \in \RR^{d} ,
\end{array}\right.
\end{equation}
where $\eps$ is a positive parameter and $\theta(\eps) \underset{\eps \to 0}{\longrightarrow} 0$ to be chosen later, and $Q$ is the collisional linear Boltzmann operator given by
\begin{equation}
Q(f) := \int_{\RR^{d}} [\sigma(x,v,v')f(v') - \sigma(x,v',v)f(v)] \ \ud v'
\end{equation}
with a non-negative \emph{collision kernel} $\sigma =\sigma(x,v,v') \geqslant 0$. The drift term $j^\eps_F$ is given by the first moment, when it exists, of the equilibrium of $Q$. It will be defined later in the statement of the main result.

For non-negative initial data $f_0$,  the unknown $f^\eps(t,x,v) \geqslant 0$ can be interpreted as the density of particles occupying, at time $t \geqslant 0$, the position $x \in \RR^d$ with velocity $v \in \RR^d$,   since for the following rescaling
$$ t'=\frac{t}{\theta(\eps)} \quad \mbox{ and } \quad x'=\frac{x}{\eps} , $$
the function $f^\eps$ can be seen as solution of the family of equations (without primes)
\begin{equation}\label{equation de BL}
\left\{\begin{array}{l}
 \partial_t f + (v-j^\eps_F)\cdot \nabla_x f = Q_\eps(f), \quad   x \in \RR^d , \ v \in \RR^d,\ t \geqslant 0 ,\\
\\
    f(0,x,v) = f_0(x,v), \hspace{2.08cm} x \in \RR^d , \ v \in \RR^d ,
\end{array}\right.
\end{equation}
where  $$ Q_\eps(f) := \int_{\RR^{d}} [\sigma(\eps x,v,v')f(v') - \sigma(\eps x,v',v)f(v)] \ \ud v' .  $$
Note that for $\eps$ small enough, the time and position variables considered in equation \eqref{BL-eps} are much larger than those in equation \eqref{equation de BL}. Thus, the study of the behavior of $f^\eps$ when $\eps$ tends to $0$ can be seen as a kind of hydrodynamic limit for the kinetic equation \eqref{equation de BL}, whose parameter $\eps$ designates the mean free path, and plays the role of Knudsen's number here.   \medskip

\ni \textbf{Some classical notations.} The collision operator $Q$ can be decomposed into a ``gain'' term and a ``loss'' term as follows:
$$ Q(f) = Q^+(f) - Q^-(f)  ,$$
with
$$ Q^+(f) := \int_{\RR^{d}} \sigma(x,v,v')f(v')\ud v'  \qquad \mbox{and} \qquad Q^-(f) := \nu(x,v) f  ,$$
where $\nu$ is the \emph{collision frequency} defined by
$$
\nu(x,v) := \int_{\RR^{d}} \sigma(x,v',v)\ud v' = \frac{Q^+(F)}{F(v)} ,
$$
and $F$ is the equilibrium of $Q$, a positive function satisfying
\begin{equation}\label{equilibre F}
Q(F)=0 \quad \mbox{ and } \quad \int_{\RR^d} F(v)\ \mathrm{d}v = 1 . 
\end{equation}
\begin{remark}
\item \begin{enumerate}
\item This splitting, $Q = Q^+ - Q^-$, is not possible in general because the integrals defining each piece diverge, but this is not the case in this paper.  
\item Such an equilibrium $F$ satisfying \eqref{equilibre F} exists under certain assumptions on the collision kernel $\sigma$, thanks to the Krein-Rutman theorem  (see \cite{DeGoPo} for more details).  A particular case in which this condition is satisfied is when $\sigma$ is such that
$$
\forall x, v, v' \in  \RR^d , \quad \sigma(x,v,v') = \mathrm{b}(x,v,v')F(v) \quad \mathrm{with}  \quad \mathrm{b}(x,v',v) = \mathrm{b}(x,v,v'), 
$$
for some non-negative $\mathrm{b} \in  L^1_{\mathrm{loc}}(\RR^{3d})$.  
\item We return to the assumptions on $\sigma$ and $F$ with more details in the next subsection. See Remarks \& Examples 3.1 in \cite{MMM} for a discussion of assumptions about the cross section, and Lemma 6.1 in the same reference for the proof of statements in some cases. 
\end{enumerate}
\end{remark}

Formally, passing to the limit when $\eps \to 0$ in equation \eqref{BL-eps}, we obtain that the limit $f^0$ is in the kernel of $ Q$ which is spanned by the equilibrium $ F $, which means that $f^0=\rho(t,x)F(v)$. Thus, it amounts to find the equation satisfied by the density $\rho$.    This question of approximation of kinetic equations by macroscopic equations has a long history, dating back to the pioneering works of E. Wigner \cite{wigner1961nuclear}, A. Bensoussan et al. \cite{BeLi}, as well as E.W. Larsen and J.B. Keller \cite{LaKe}. Since then, numerous papers have addressed this topic (for further references, see the papers by C. Bardos et al. \cite{BaSaSe}, and P. Degond et al. \cite{DeGoPo}). The resulting equations fall into two categories, depending on the rate of decrease in velocity of the equilibrium $F$.  It is well-known (see \cite{DeGoPo} for instance) that when $F$ decreases \emph{rapidly} for large velocities (such as when $F$ follows a Maxwellian distribution function), the distribution function $f^\eps$ converges to $\rho(t,x)F(v)$ as $\eps$ goes to $0$, for the classical scaling $\theta(\eps) = \eps^2$,  with $\rho$ being the solution of the diffusion equation
\begin{equation}\label{eq de diff BL}
 \pa_t \rho - \na_x \cdot (D \na_x \rho) = 0 ,
\end{equation}
where 
$$  D = \int_{\RR^d} (v \otimes v) \frac{F(v)}{\nu(v)} \ud v .$$
When $F$ decreases \emph{slowly} and it is a  heavy tail distribution function, satisfying
\begin{equation}\label{F(v) sim |v|^(-alpha-d) intro BL}
F(v) \sim \frac{\kappa}{|v|^{d+\alpha}} \quad \mbox{ as } \quad |v| \to \infty 
\end{equation}
for some $\alpha > 0$, the previous diffusion matrix $D$, given by \eqref{eq de diff BL} might be infinite.  In that case, the diffusion limit leading to \eqref{eq de diff BL} breaks down,  which means that the choice of time scale $\theta(\eps) = \eps^2$ was inappropriate.  It has been shown in \cite{MMM} and \cite{M},  using different methods and under various assumptions, that in such cases, the appropriate time scale involves the parameter $\alpha$, which appears in the equilibrium. More specifically, for $\theta(\eps) = \eps^\gamma$ with $\gamma$ depends on $\alpha$,  the following fractional diffusion equation is obtained
\begin{equation}\label{eq diff frac BL}
\pa_t \rho + \kappa (-\Delta_x)^{\frac{\gamma}{2}} \rho = 0 ,
\end{equation}
where the fractional Laplacian appearing in the previous equation can be defined by the following singular integral
$$ (-\Delta_x)^{\frac{\gamma}{2}} u(x) := c_{\alpha,d} \ \mathrm{PV} \int_{\RR^d} \frac{u(x)-u(y)}{|x-y|^{d+\gamma}} \ \ud y .$$ 
Let us now discuss in more detail the contexts of the last two references cited.  In \cite{MMM},  the authors addressed the problem in the space homogeneous case (that is with $\sigma$  independent of $x$).  They prouved that when $F$ satisfies \eqref{F(v) sim |v|^(-alpha-d) intro BL} and the collision frequency $\nu$ satisfies
$$ \nu(v) \sim  \nu_0 |v|^{\beta} \quad \mbox{ as } \quad |v| \to \infty , $$
then for $\alpha > 0$ and $\beta < \min(\alpha;2-\alpha)$,  and by taking $\gamma:=\frac{\alpha-\beta}{1-\beta}$,  the scaling  $\theta(\eps)=\eps^\gamma$ leads to the previous fractional diffusion equation.  While in \cite{M}, the same problem has been addressed but for a cross section which depends on the position variable $x$ but assuming that the collision frequency is bounded, more precisely
\begin{equation}\label{sigma entre nu_i F BL}
 0<\nu_1 F(v) \leqslant \sigma(x,v,v') \leqslant \nu_2 F(v),
\end{equation}
i.e.  for $\beta=0$ compared to the case of \cite{MMM}, which gives in particular,
\begin{equation}\label{nu entre nu_i}
0< \nu_1  \leqslant \nu(x,v) \leqslant \nu_2 .
\end{equation}
With the additional assumption
\begin{equation}\label{nu et nu_0 intro BL}
\nu(x,v) \sim \nu_0(x) \quad \mbox{ as } \ |v| \to \infty \quad \mbox{ with } \quad \nu_1 < \nu_0(x) < \nu_2 ,
\end{equation}
the same type of macroscopic equation as \eqref{eq diff frac BL} was found, with an elliptic operator $\mathcal{L}$ of the same order as $(-\Delta_x)^{\alpha/2}$
$$ \mathcal{L}(\rho) := \mathrm{PV} \int_{\RR^{d}}\eta(x,y)\frac{\rho (x) - \rho(y)} {\vert x - y \vert ^{d + \alpha}} \ \ud y , $$
with
$$ \eta(x,y) = \nu_0(x)\nu_0(y) \int_0^\infty z^\alpha e^{ -z\int_0^1\nu_0(sx+(1-s)y) \ud s} \ \ud z .$$

 Towards the end of this section, we mention some other references to previous results in the context of anomalous and fractional diffusion. The earliest work dates back to \cite{borgers1992diffusion}, where the authors showed the diffusion limit of free molecular flow in thin plane channels. There are also the works of Dogb\'e and Golse on anomalous diffusion for the Knudsen gas, see for instance \cite{golse1998anomalous, dogbe1998diffusion, dogbe2000anomalous}. The term ``anomalous'' here stems from the scaling $\theta(\eps)=\eps^2 \ln(\eps^{-1})$ considered in the critical case (i.e., the limiting case for the integrability of $\int (v \otimes v)\frac{F(v)}{\nu(v)}  \mathrm{d}v$).   Regarding fractional diffusion, the first result was obtained simultaneously by A. Mellet et al. \cite{MMM}, and M. Jara et al. \cite{MKO}, using two completely different approaches (using a probabilistic approach for the latter refrence).  This was followed by the work of A. Mellet \cite{M} for heavy-tail equilibrium with cross-section satisfying Assumptions \eqref{sigma entre nu_i F BL}, \eqref{nu entre nu_i}, and \eqref{nu et nu_0 intro BL}, N. Ben Abdellah et al. \cite{BMP-AD} for a collision frequency that degenerates for small velocities with Maxwellian equilibrium, and in \cite{BMP-FD} for a Sobolev initial data to obtain strong convergence of $f^\eps$ to $\rho(t,x)F(v)$ (via a completely different method from the previous ones, based on a Hilbert expansion). Very recently, the results of \cite{MMM} have been recovered by E. Bouin and C. Mouhot \cite{BM} using a quasi-spectral method, which is also unified for the Fokker-Planck and L\'evy-Fokker-Planck models.  

 To the author's knowledge, all existing results on the diffusion limit of the linear Boltzmann equation have been established under the assumptions of a symmetric collision frequency and equilibrium. In this paper, no symmetry is assumed for these two functions, and as a consequence, a drift term must necessarily be added to the kinetic equation to ensure the convergence of $f^\eps$. However, in the case of the Fokker-Planck equation, the diffusion limit for an asymmetric equilibrium has been obtained very recently. This was achieved through a purely probabilistic approach in \cite{FT} for the case $d \geqslant 2$, and via a spectral approach in \cite{dechicha2024spectral} in any dimension $d \geqslant 1$. This last paper extends the results of \cite{dechicha2023construction, DP-2}, which generalize \cite{LebPu}.

 The goal of this paper is to show that a similar phenomenon, such as the fractional diffusion limit, can arise, and that a macroscopic equation with an elliptic operator equivalent to the fractional Laplacian is obtained from \eqref{BL-eps} in situations in which the equilibrium is heavy-tailed (satisfies \eqref{F(v) sim |v|^(-alpha-d) intro BL} for some $\alpha>0$),  the cross-section $\sigma$ depends on the spatial variable $x$, and the collision frequency $\nu$ degenerates for large velocities
$$ \nu(x,v) \sim |v|^{\beta}\nu_0(x) \quad \mbox{ as } \quad |v| \to \infty \quad \mbox{ with } \quad \nu_1 < \nu_0(x) < \nu_2 , $$
for some $\beta \in \RR$ (with $\beta < \min(\alpha;2-\alpha)$), without assuming symmetry, either for the equilibrium $F$ or for the collision frequency $\nu$.  For this, we will adapt the method used in \cite{M}, slightly modifying the auxiliary problem in order to handle the situation when $\nu$ degenerates for $|v|$ large enough and considering for $\varphi$ smooth,  the equation 
$$ |v|^{-\beta} \nu(x,v)\chi^\eps - \eps |v|^{-\beta}  v \cdot \nabla_x \chi^\eps = |v|^{-\beta}  \nu(x,v) \varphi(t,x) $$ 
for $|v| \geqslant 1$ and keeping the same problem for $|v| \leqslant 1$. The function $\chi^\eps(t,x,v)$, which will converge to $\varphi$ when $\eps$ goes to $0$, will be considered as the test function in the weak formulation of equation \eqref{BL-eps}.  Also, due to the degeneracy of $\nu$, the operator $Q$ has no spectral gap and the passing to the limit in most of terms in the weak formulation is treated differently to \cite{M}.  In addition, we have more terms to handle in this paper due to the presence of the drift term.

\subsection{Main result}\label{main thm}
First,  we denote by $L^p_{\omega}$ the space $L^p$ endowed with the measure $\omega \ud v$ (or $\omega \ud v \ud x$). For example, $L^2_{\nu F^{-1}}(\RR^d)$ is the space defined by
$$ L^2_{\nu F^{-1}}(\RR^d) := \left\{ f : \RR^d \longrightarrow \RR ; \ \int_{\RR^d} |f|^2 \frac{\nu}{F} \ \ud v < \infty \right\} . $$
Before stating the result, we need to make the conditions on $\sigma$ and $F$ precise.  
\begin{assump}
The first two assumptions are very standard:
\item (A1) The cross section $\sigma(x,v,v')$ is non negative and locally integrable on $\RR^{2d}$ for all $x$. 
The collision frequency $ \nu(x,v) = \int_{\RR^{d}} \sigma(x,v',v)\ud v'$ satisfies 
$$ \nu(x,v) > 0 , \qquad \forall x,v \in \RR^{d} .$$
\item (A2) There is a function $F(v) \in L^1_\nu(\RR^{d})$ independent of $x$ such that
$$ Q(F) = 0  .$$
Then,  $ \ Q^+(F) = \nu F$.  Furthermore, the function $F$ is positive and normalized to $1$:
 $$ F(v) > 0 \quad \mathit{for}\ \mathit{ all }\  v \in \RR^{d} \quad \mathit{ and } \quad \int_{\RR^{d}} F(v) \ud v =  1  .  $$
 \end{assump}
 \begin{remark}
Note that we do not assume symmetry in this paper, unlike \cite{MMM} and \cite{M}, where they assume that $F(-v) = F(v)$ and $\nu(-v) = \nu(v)$ for all $v$ in $\RR^d$.
 \end{remark}

 The next assumptions concern the behavior of $F$ and $\nu$ for large $|v|$. We combined the assumptions of \cite{MMM} and \cite{M} to work with a more general collision frequency.
\begin{assump}
\item (B1) There exists  $\alpha > 0$ and a constant  $\kappa > 0$ such that
$$ |v|^{\alpha + d}F(v) \longrightarrow \kappa \quad \mathit{as} \quad |v| \rightarrow \infty .$$
\item (B2) There exists $\nu_1$ and $\nu_2$ two positive constants,  a function $\nu_0(x)$ and a constant $\beta \in  \mathbb{R}$ such that
$$ \nu_1 \langle v \rangle^\beta \leqslant \nu(x,v) \leqslant \nu_2 \langle v \rangle^\beta , \qquad \forall x,v \in \RR^{d}  $$
\emph{and }
$$ |v|^{- \beta}  \nu (x,v) \longrightarrow \nu_0(x)\quad \mathit{as} \quad |v| \longrightarrow \infty ,  \quad \mathit{ uniformly}\ \mathit{ with }\ \mathit{ respect }\ \mathit{ to }\ x.   $$ 
We assume that $\nu$ is $C^1$ with respect to $x$ and
   $$ \big\| \langle v \rangle^{-\beta} \partial_x \nu (x,v) \big\|_{L^{\infty}(\RR^{2d})} \ \leqslant C   .$$
\item \emph{(B3) Finally, we assume that there exists a positive constant $M$ such that}
   $$ \int_{\RR^{d}} F'\frac{\nu}{\mathrm{b}} \ \ud v' +  \left( \int_{\RR^{d}}\frac{F'}{\nu'}\frac{\mathrm{b}^2}{\nu^2} \ \ud v' \right) ^{\frac{1}{2}} \leqslant M , \qquad \forall x,v \in \RR^{d} ,$$
 with $\mathrm{b} := \mathrm{b}(x,v,v') := \sigma(x,v,v')F^{-1}(v)$, $F'=F(v')$ and $\nu':=\nu(x,v')$. 
\end{assump}  
\begin{remark}
This last assumption, (B3), concerns the collision operator $Q$, so that it is bounded and has the property of weighted coercivity.
\end{remark}

\begin{theorem}[Fractional diffusion limit for the linear Boltzmann equation]\label{main BL}
\item Assume (A1-A2) and (B1-B2-B3) with $\alpha > 0$ and $\beta < \min\{\alpha;2-\alpha\}$.  We define
 $ \gamma := \frac{\alpha - \beta}{1 - \beta}$, and the drift term $j^\eps_F$ by
 \begin{equation}\label{j-eps}
j^\eps_F :=\left\{\begin{array}{l}
 0 \hspace{3.26cm} \mathrm{if} \ \alpha < 1  , \\
\\
    \ds \int_{|v| \leqslant \eps^{-\frac{1}{1-\beta}}} vF\ud v \qquad \mathrm{if} \ \alpha = 1, \\
    \\
    \ds \int_{\RR^d} vF\ud v \hspace{1.89cm} \mathrm{if} \  \alpha > 1 .
\end{array}\right.
\end{equation}
Let $f^\varepsilon(t,x,v)$ be the solution of \eqref{BL-eps} with $\theta(\varepsilon)=\varepsilon^\gamma$ and non-negative $f_0 \in L^2_{F ^{-1}}(\RR^{2d})$.  Then, $f^\varepsilon$ converges weakly star in $L^\infty(0,T;L^2_{\nu F^{-1}}(\RR^{2d}))$ to $\rho(x ,t)F(v)$ with $\rho$ solution to
   \begin{equation}
   \left\{\begin{array}{l}
 \partial_t \rho + \kappa \mathcal{L}(\rho) = 0 ,  \\
 \\
 \ds  \rho(x,0) = \rho_0(x) := \int_{\RR^d} f_0(x,v) \ \ud v ,
\end{array}\right.
\end{equation}    
where $\mathcal{L}$ is an elliptic operator of order $\gamma$ defined by
\begin{equation}\label{definition de L BL}
\mathcal{L}(\rho) := \frac{1}{1-\beta} \ \mathrm{PV} \int_{\RR^{d}}\eta(x,y)\frac{\rho (x) - \rho(y)}{\vert x - y \vert ^{d + \gamma}} \ \ud y , 
\end{equation}
with 
$$ \eta(x,y) := \nu_0(x)\nu_0(y) \int_0^\infty z^\gamma e^{ -z\int_0^1\nu_0(sx+(1-s)y) \ud s} \ \ud z .$$
\end{theorem}
\begin{remark}
\item \begin{enumerate}
\item  Note that $\beta < 1$ and $ \gamma < 2$ for $\alpha > 0$ and $\beta < \min\{\alpha;2-\alpha\}$. 
\item From $(B2)$,  there exists $\eta_1$ and $\eta_2$ such that
$$ 0 < \eta_1 \leqslant \eta(x,y) \leqslant \eta_2 < \infty .$$ 
In particular, the operator $\mathcal{L}$ has the same order as the fractional Laplacian operator $(-\Delta_x)^{\gamma/2}$. Moreover, for $\sigma$ independent of $x$, the function $\eta$ is constant and formula \eqref{definition de L BL} gives exactly $(-\Delta_x)^{\gamma/2}$ up to a multiplicative positive constant.
\item The operator $\mathcal{L}$ is self-adjoint since $\eta(x,y) = \eta(y,x)$.
\end{enumerate}    
\end{remark}

\begin{remark} In the particular case $F(-v)=F(v)$ and $\nu(x,-v)=\nu(x,v)$ for all $v \in \RR^d$, and for a cross section $\sigma$ independent of $x$, we recover Theorem 3.2 in \cite{MMM}, and for $\beta = 0$ we recover the results of \cite{M}.
\end{remark}
\section{Preliminaries results}
We start by recalling some classical properties of the collision operator $Q$.
\subsection{Properties of the operator $Q$}
Under Assumptions (A1-A2) and (B3), we have the following proposition.
\begin{proposition}\label{coercivite BL}
Let $f$ and $g$ two functions in $L^2_{\nu F^{-1}}(\RR^d)$.  Then we have the following assertions:
\begin{enumerate}
\item The operator $Q:L^1_\nu \longrightarrow L^1$ is bounded and conservative, so equation \eqref{equation de BL} preserves the total mass of the distribution$f$ and one has
$$\int_{\RR^d} Q(f)\ \ud v= 0,\qquad \forall   f \in L^1_{\nu}(\RR^d).$$
\item The operator $\frac{1}{\nu}Q$ is bounded in $L^2_{\nu F^{-1}}$ and $Q$ is dissipative.   Moreover, 
\begin{equation}\label{Poincare BL}
\int_{\RR^{d}} Q(f)f \ \frac{\ud v}{F}  \leqslant -\frac{1}{2M}\int_{\RR^{d}} \vert f - \rho F \vert ^2 \ \frac{\nu \ud v}{F}  , \qquad \forall f \in L^2_{\nu F^{-1}} .
\end{equation} 
where  $\ds \rho := \int_{\RR^{d}} f\ \ud v$ and $M$ is the constant given in Assumption (B3).
 \end{enumerate}
 \end{proposition}
 
 \begin{remark}
 \item\begin{enumerate}
 \item It should be noted that for an unbounded collision frequency, which degenerates for large velocities, for example, the integrals in inequality \eqref{Poincare BL} are conducted with two different measures, namely $\frac{\ud v}{F}$ and $\frac{\nu \ud v}{F}$. Thus, the operator $Q$ does not have a spectral gap for the type of cross-sections or collision frequencies considered in this work.
 \item Note that for $F \in L^2_{\nu F^{-1}}(\RR^d)$,  we get the following inclusion $\ L^2_{\nu F^{-1}}(\RR^d) \subset L^1_\nu(\RR^d)$. 
 \end{enumerate}
\end{remark}
 
\ni \bpp \ref{coercivite BL}.   We adapt the proof of Lemma 4.1 in \cite{MMM},  whose authors had adapted the proof of Propositions 1 \& 2 from \cite{DeGoPo}.  Recall that $Q(f) = Q^+(f) - \nu f$. Thus, to show that $\frac{1}{\nu} Q$ is bounded in $L^2_{\nu F^{-1}}$, it amounts to showing that $\frac{1}{\nu } Q^+$ is bounded.   \\
Let $f \in L^2_{\nu F^{-1}}$.  By using the fact that $Q^+(F)=\nu F$ and $\big(\int \sigma f \ud v\big)^2 \leqslant \big(\int \sigma \frac{f^2}{F} \ud v\big)\big(\int \sigma F \ud v\big)$, we write
\begin{align*}
\left\| \frac{1}{\nu}Q^+(f) \right\| ^2_{L^2(\nu F^{-1})} &= \int_{\RR^{d}}\frac{\big| Q^+(f)\big|^2}{\nu F} \ud v \\
&\leqslant \int_{\RR^{d}}\frac{1}{\nu F} \bigg(\int_{\RR^{d}}\sigma(v,v')F(v')\ud v' \bigg)\bigg(\int_{\RR^{d}}\sigma(v,v')\frac{\vert f(v') \vert ^2}{F(v')}\ud v' \bigg)\ud v\\
& =  \int_{\RR^{d}}\frac{1}{\nu F} Q^+(F)\int_{\RR^{d}}\sigma(v,v')\frac{\vert f(v') \vert ^2}{F(v')}\ud v' \ud v  = \iint_{\RR^{2d}} \sigma(v,v')\frac{\vert f(v') \vert ^2}{F(v')}\ud v' \ud v \\
& = \int_{\RR^{d}} \nu(v')\frac{\vert f(v') \vert ^2}{F(v')}\ud v' = \| f \|^2_{L^2_{\nu F^{-1}}} .
\end{align*}
Now, let's prouve  inequality \eqref{Poincare BL}.  We have
 $$ \int_{\RR^{d}} Q(f)\frac{f}{F} \ud v = \iint_{\RR^{2d}} \sigma(v,v')f'\frac{f}{F}\ud v' \ud v - \int_{\RR^{d}}\nu \frac{f^2}{F}\ud v . $$
On the one hand, 
$$ \int_{\RR^{d}}\nu(v) \frac{f^2}{F}\ud v = \iint_{\RR^{2d}} \sigma(v',v)F\frac{f^2}{F^2}\ud v \ud v' = \iint_{\RR^{2d}} \sigma(v,v')F' \frac{f'^2}{F'^2}\ud v \ud v' .  $$
On the other hand,
 $$ \int_{\RR^{d}}\nu(v) \frac{f^2}{F}\ud v = \int_{\RR^{d}} Q^+(F)\frac{f^2}{F^2}\ud v = \iint_{\RR^{2d}} \sigma(v,v')F' \frac{f^2}{F^2} .$$ 
Hence, 
 $$ \int_{\RR^{d}} Q(f)\frac{f}{F} \ud v = -\frac{1}{2} \iint_{\RR^{2d}} \sigma(v,v')F'\left[ \frac{f}{F} - \frac{f'}{F'}\right] ^2\ud v' \ud v .  $$ 
Moreover, we have
 $$ fF' - f'F = \left( \frac{f}{F} - \frac{f'}{F'} \right) FF' .$$ 
Then,
 $$ g := f - \rho F =  \int_{\RR^{d}}(fF' - f'F)\ud v' = \int_{\RR^{d}} \left( \frac{f}{F} - \frac{f'}{F'} \right) FF'\ud v' ,$$
which implies that
 $$ g^2 \leqslant \left( \int_{\RR^{d}} \sigma(v,v') F' \left( \frac{f}{F} - \frac{f'}{F'} \right)^2 \ud v'\right) \left( \int_{\RR^{d}}\frac{F^2}{\sigma(v,v')} F' \ud v'\right)  .  $$
Thus, by integrating this last inequality against $\frac{\nu}{F}$ and using Assumption (B3), we obtain
 $$ \int_{\RR^{d}} g^2 \frac{\nu}{F} \ \ud v \leqslant \left( \sup_{v \in \RR^{d}} \left( \nu \int_{\RR^{d}}\frac{F}{\sigma} F' \ud v'\right)\right) \left( \iint_{\RR^{2d}}\sigma(v,v') F' \left( \frac{f}{F} - \frac{f'}{F'} \right)^2 \ud v'\ud v\right) .$$
Therefore, 
 $$ \int_{\RR^{d}} g^2\frac{\nu}{F} \ud v \leqslant -2M \int_{\RR^{d}} Q(f)\frac{f}{F} \ud v .$$
Hence the inequality of the second item of Proposition \ref{coercivite BL} holds true.
\epp

\begin{remark} If we had considered the integral $\int_{\RR^{d}} g^2 \frac{ \ud v}{F}$ instead of $\int_{\RR^{d}} g^2 \frac{\nu}{F} \ \ud v$, we would have obtained $\sup_{v \in \RR^{d}} \left(\int_{\RR^{d}}\frac{F F' }{\sigma} \ud v'\right)$ in the right-hand side of the penultimate inequality, which may not be bounded for typical physical examples where $\nu(v) \sim |v|^\beta$.
\end{remark}
\subsection{Compactness lemma}
As a consequence of the previous proposition, we get the following compactness lemma.
\begin{lemma}\label{lemme de compacite BL}
For initial datum $f_0 \in L^2_{F^{-1}}(\RR^{2d})$ and a positive time $T$,
\begin{enumerate}
\item The solution $f^\eps$ of \eqref{BL-eps} is bounded in $L^{\infty}\big(0,T;L^2_{F^{-1}}(\RR^{2d})\big)$ uniformly with respet to $\eps$.  Moreover, 
\begin{equation}\label{g^eps < theta(eps) f_0}
 M \| f_\eps \|^2 _{L^2_{F^{-1}}} \theta(\eps) + \| f^\eps - \rho^\eps F \|^2 _{L^{\infty}(0,T,L^2_{\nu F^{-1}})} \leqslant M \| f_0 \|^2 _{L^2_{F^{-1}}} \theta(\eps) . 
\end{equation}
\item The density  $\ds \rho^{\varepsilon} := \int_{\RR^{d}} f^\eps \ \ud v$ is bounded in $L^{\infty}\big(0,T;L^2(\RR^{d})\big)$ uniformly with respet to $\eps$ and one has
\begin{equation}\label{rho^eps < f_0}
 \| \rho^{\varepsilon} \| _{L^{\infty}(0,T;L^2)} \leqslant \| f_0 \| _{L^2_{F^{-1}}} . 
\end{equation}
\end{enumerate}
In particular,  $\rho^\eps$ converges weakly star in $L^{\infty}(0,T;L^2(\RR^d))$ to $\rho$ and $f^\eps$ converges weakly star in $ L^{\infty}(0,T; L^2_{\nu F^{-1}}(\RR^{2d}))$ to $f = \rho (x,t) F(v)$.
\end{lemma}

\ni \bp 1.  By integrating equation \eqref{BL-eps} multiplied by $f^\eps/F$ and using inequality \eqref{coercivite BL}, we obtain
\begin{align*}
\frac{1}{2}\frac{\ud}{\ud t}\int_{\RR^{2d}}  \frac{|f^{\varepsilon}|^2}{F} \ud x \ud v &=\frac{1}{\theta(\varepsilon)}\int_{\RR^{2d}} Q(f^{\varepsilon}) \frac{ f^{\varepsilon}}{F} \ud x \ud v \\
&\leqslant -\frac{1}{2M\theta(\varepsilon)} \int_{\RR^{2d}} \big| f^{\varepsilon} - \rho^{\varepsilon} F \big|^2 \frac{\nu}{F} \ud x \ud v .
\end{align*}
Therefore,
$$ \frac{1}{2}\int_{\RR^{2d}} |f^\eps(t,x,v)|^2 \ \frac{\ud x \ud v}{F}  + \frac{1}{2M\theta(\varepsilon)} \int_0^t\int_{\RR^{2d}} \big| f^\eps - \rho^\eps F \big|^2 \frac{\nu}{F} \ud x \ud v  \leqslant \frac{1}{2}\int_{\RR^{2d}} | f_0(x,v)|^2 \ \frac{\ud x \ud v}{F} .  $$ 

This last inequality shows that $f^\eps$ is uniformly bounded in $L^{\infty}(0,T,L^2_{F^{-1}})$ and gives also  inequality \eqref{g^eps < theta(eps) f_0}.   \smallskip

\ni 2.  For the second point,  we write
$$ \int_{\RR^{d}} |\rho^\eps |^2 \ud x = \int_{\RR^{d}} \bigg| \int_{\RR^{d}} f^{\varepsilon} \ud v \bigg| ^2 \ud x \leqslant \bigg(\int_{\RR^{d}} F \ud v\bigg)\bigg(\int_{\RR^{2d}} | f^{\varepsilon} |^2 \frac{\ud v \ud x}{F}\bigg)  = \| f^\eps \|^2_{L^2_{F^{-1}}} .$$
\ep

\subsection{The auxiliary problem}
We adapt the method used in \cite{M} and this step is the key of the proof. It consists in correcting the error on the test function in the weak formulation.  Indeed,  for a test function $\varphi \in \mathcal{D}([0,+\infty)\times\RR^{d})$, we introduce $\chi^{\varepsilon}(t,x,v)$ solution to
\begin{equation}\label{eq de chi^eps}
\lceil v \rfloor^{-\beta} \nu(x,v)\chi^\eps - \eps \lceil v \rfloor^{-\beta} v \cdot \nabla_x \chi^\eps = \lceil v \rfloor^{-\beta} \nu(x,v) \varphi(t,x) ,
\end{equation}
where 
$$ \lceil v \rfloor^{-\beta} :=   \left\{\begin{array}{l}
1  \hspace{1.04cm} \mbox{for} \quad |v| \leqslant 1 , \\
  |v|^{-\beta}\quad \mbox{for} \quad  |v| \geqslant 1 .
\end{array}\right.$$
Note that for $\sigma > 0$,
$$ \lceil v \rfloor^{\sigma} \leqslant \langle v \rangle^{\sigma} \leqslant 2^{\sigma/2} \lceil v \rfloor^{\sigma} .$$
The trick which makes it possible to generalize the work given in \cite{M} to a cross section which degenerates for large velocities is in the term $\lceil v \rfloor^{-\beta}$ introduced in the previous equation. It allows to renormalize the collision frequency $\nu$ in some sense, and it is as if we consider particles moving with a velocity $\lceil v \rfloor^{-\beta} v$.  \medskip

We can integrate equation \eqref{eq de chi^eps} and its solution $\chi^\eps$ is given explicitly, according to the function $\varphi$, by the following formula
$$
\chi^{\varepsilon}(t,x,v) = \int_0^\infty \lceil v \rfloor^{-\beta} \nu(x+\varepsilon\lceil v \rfloor^{-\beta} vz,v) \ e^{-\int_0^z \lceil v \rfloor^{-\beta} \nu(x+\varepsilon\lceil v \rfloor^{-\beta} vs,v) \ud s} \varphi(t,x+\varepsilon\lceil v \rfloor^{-\beta} vz) \ \ud z .
$$
In order to simplify the writing and to avoid long expressions, we denote by $\tilde \nu$ and $\tilde \varphi$ the two following functions
$$ \tilde \nu_\eps(x,v,z) :=  \nu(x+\varepsilon\lceil v \rfloor^{-\beta} vz,v)  \qquad \mbox{and} \qquad \tilde \varphi_\eps(t,x,v,z) := \varphi(t,x+\varepsilon\lceil v \rfloor^{-\beta} vz) .$$
Then,
\begin{equation}\label{formule chi^eps}
\chi^{\varepsilon}(t,x,v) = \int_0^\infty \lceil v \rfloor^{-\beta} \tilde \nu_\eps(x,v,z) \ e^{- \int_0^z \lceil v \rfloor^{-\beta} \tilde \nu_\eps(x,v,s) \ud s} \tilde \varphi_\eps(t,x,v,z) \ \ud z .
\end{equation}
The function $\chi^\eps$ is smooth. Furthermore,  thanks to the identity

$$ \int_0^\infty \lceil v \rfloor^{-\beta} \tilde \nu_\eps(x,v,z) \ e^{- \int_0^z \lceil v \rfloor^{-\beta} \tilde \nu_\eps(x,v,s) \ud s}  \ \ud z  = \int_0^\infty e^{-u} \ \ud u = 1 , $$
and Assumption (B2), we write
\begin{align*}
|\chi^{\varepsilon} - \varphi| &= \bigg |\int_0^\infty \lceil v \rfloor^{-\beta} \tilde \nu_\eps(x,v,z) \ e^{- \int_0^z \lceil v \rfloor^{-\beta} \tilde \nu_\eps(x,v,s) \ud s} [\tilde \varphi_\eps(t,x,v,z)-\varphi(t,x)]  \ud z \bigg| \\
&\leqslant \nu_2 \int_0^\infty e^{-\nu_1z}\varepsilon |v|\lceil v \rfloor^{-\beta}  z \| \na_x\varphi\|_{L^{\infty}} \ud z = \frac{\nu_2}{\nu_1} \eps |v| \lceil v \rfloor^{-\beta} \| \na_x \varphi\|_{L^{\infty}} .
\end{align*}
Thus,  $\chi^\eps$ is bounded in $L^{\infty}_{t,x}$ and converges to $\varphi$,  uniformly with respect to $x$ and $t$.  However, this convergence is not uniform with respect to $v$.  Moreover, we have the following lemma.
\begin{lemma}\label{chi et phi}
Let $\varphi \in \mathcal{D}([0,\infty)\times\RR^{d})$ and let $\chi^{\varepsilon}$ the function given by \eqref{eq de chi^eps}.  Then,
\begin{equation}\label{F(chi-phi)-->0}
\int_0^\infty \int_{\RR^{2d}} F(v)\big[\chi^\eps(t,x,v) - \varphi(t,x) \big]^2 \ud x \ud v \ud t \underset{\eps \to 0}{\longrightarrow} 0
\end{equation}
and
\begin{equation}\label{F(chi_t - phi_t)-->0}
\int_0^\infty \int_{\RR^{2d}}  F(v)\big[\pa_t\chi^\eps(t,x,v) - \pa_t\varphi(t,x) \big]^2 \ud x \ud v \ud t \underset{\eps \to 0}{\longrightarrow} 0 .
\end{equation}
Moreover, there is a constant $C>0$ such that
\begin{equation}\label{chi bornee par phi}
\| \chi^\eps \|_{L^2_{F}((0,\infty)\times\RR^{2d})} \leqslant C \| \varphi\|_{L^2((0,\infty)\times\RR^{d})}  
\end{equation}
and 
\begin{equation}\label{chi_t bornee par phi_t}
\| \pa_t\chi^\eps \|_{L^2_{F}((0,\infty)\times\RR^{2d})} \leqslant C \| \pa_t\varphi\|_{L^2((0,\infty)\times\RR^{d})} .
\end{equation}
\end{lemma}

\ni This lemma is in a way the equivalent of Lemmas 3.1 and 4.2 in \cite{M} but not quite the same thing, since in \cite{M} we have
$$ \int_{\RR^d} F(v)\big[\chi^\eps(t,x,v) - \varphi(t,x) \big] \ud v \underset{\eps \to 0}{\longrightarrow} 0 $$
and the same for the derivative.  The proof of Lemma \ref{chi et phi} is different from the one given in \cite{M}, since in the formula of $\chi^\eps$ involves $\lceil v \rfloor^{-\beta}$, and even more so in its derivative.  We obtain only weak limits instead of the bounds by $\eps$.  \\

\ni \bpl \ref{chi et phi}.  We are going to prove limit \eqref{F(chi_t - phi_t)-->0}, and \eqref{F(chi-phi)-->0} is done exactly the same way.  We have
$$ \| \partial_t\chi^{\varepsilon} - \partial_t\varphi \|_{L^2_F}^2 = \int_0^{\infty}\int_{\RR^{2d}} F\big[ \partial_t\chi^{\varepsilon} - \partial_t\varphi \big]^2 \ud x\ud v\ud t = \int_{\RR^{d}} F\int_0^{\infty}\int_{\RR^{d}}\big[ \partial_t\chi^{\varepsilon} - \partial_t\varphi \big]^2 \ud x \ud t\ud v . $$
On the one hand,  using Taylor expansion and Assumption (B2), we write
\begin{align*}
 \big| \partial_t\chi^{\varepsilon} - \partial_t\varphi \big| &= \bigg| \int_0^\infty \lceil v \rfloor^{-\beta} \tilde \nu_\eps(x,v,z) \ e^{-\int_0^z \lceil v \rfloor^{-\beta} \tilde \nu_\eps(x,v,s) \ud s} \big[\partial_t \tilde \varphi_\eps(t,x,v,z) - \partial_t\varphi(t,x) \big] \ud z \bigg|  \\
 &\leqslant \int_0^{\infty}\nu_2 e^{-\nu_1 z}\varepsilon |v| \lceil v \rfloor^{-\beta} z \int_0^1 | \partial_t \na_x\varphi(t,x+\eps v \lceil v \rfloor^{-\beta}zs)| \ \ud s \ud z  .
 \end{align*}
 Then,
$$
 \big| \partial_t\chi^{\varepsilon} - \partial_t\varphi\big|^2 \leqslant \nu_2^2{\varepsilon}^2 \lceil v \rfloor^{2(1-\beta)} \bigg(\int_0^{\infty}z^2e^{-\nu_1 z}  \bigg) \bigg( \int_0^{\infty} e^{-\nu_1 z}\int_0^1 \big|\partial_t \na_x \tilde \varphi_\eps(t,x,v,zs)\big| ^2\ud s\ud z \bigg) .
$$
Therefore, 
$$  \int_0^{\infty}\int_{\RR^{d}} \big| \partial_t\chi^{\varepsilon} - \partial_t \varphi  \big|^2 \ud x \ud t \leqslant C {\varepsilon}^2 \lceil v \rfloor^{2(1-\beta)}\| \partial_t \na_x\varphi \|^2_{L^2((0,\infty)\times\RR^{d})} ,$$
where $$ C := \nu_2^2 \int_0^{\infty}e^{-\nu_1z}\ud z\int_0^{\infty}z^2e^{-\nu_1z} \ud z . $$
Hence,  for all $v \in \RR^d$
\begin{equation}\label{F int chi-phi-->0}
F(v) \int_0^{\infty}\int_{\RR^{d}} \big| \partial_t\chi^{\varepsilon} - \partial_t \varphi  \big|^2 \ud x \ud t \underset{\eps \to 0}{\longrightarrow } 0 .
\end{equation}
On the other hand,
\begin{align*}
\big| \partial_t\chi^{\varepsilon} \big|^2 &\leqslant \int_0^\infty \lceil v \rfloor^{-\beta} \tilde \nu_\eps(x,v,z) \ e^{-\int_0^z \lceil v \rfloor^{-\beta} \tilde \nu_\eps(x,v,s) \ud s} |\partial_t \tilde \varphi_\eps(t,x,v,z)|^2 \ud z \\
&\leqslant \int_0^\infty \nu_2e^{-\nu_1z} |\partial_t \tilde \varphi_\eps(t,x,v,z)|^2  \ud z .
\end{align*}
Therefore,
\begin{equation}\label{chi_t < phi_t}
\int_0^{\infty}\int_{\RR^{d}} \big| \partial_t\chi^{\varepsilon} \big|^2 \ud x \ud t \leqslant \frac{\nu_2}{\nu_1}\| \partial_t\varphi\|^2_{L^2((0,\infty)\times\RR^{d})} .
\end{equation}
Hence,   
\begin{equation}\label{F int chi-phi in L^1}
F \ \int_0^{\infty}\int_{\RR^{d}} \big| \partial_t\chi^{\varepsilon} - \partial_t \varphi  \big|^2 \ud x \ud t \leqslant  2\frac{\nu_2}{\nu_1}\| \partial_t\varphi\|^2_{L^2}F   \in L^1_v(\RR^{d}).
\end{equation}
Thus, from \eqref{F int chi-phi-->0} and \eqref{F int chi-phi in L^1}, and by Lebesgue's theorem, limit \eqref{F(chi_t - phi_t)-->0} holds true. Inequality \eqref{chi_t bornee par phi_t} follows from \eqref{chi_t < phi_t} and $ \int_{\RR^d} F \ \ud v = 1$.
\ep

\section{Proof of the main theorem}
The aim of this section is to prove the main theorem whose proof is based on the moment method.  For that, we start by giving the weak formulation, taking the solution of the auxiliary problem as a test function.  \medskip

\ni \textbf{Weak formulation.} We recall that $f^\eps$ solves \eqref{BL-eps} with $\theta(\eps) = \eps^{\gamma}$.  By multiplying equation \eqref{BL-eps} by $\chi^{\varepsilon}$ and integrating it with respect to $x$, $v$ and $t$, we obtain:
\begin{align*}
&-\varepsilon^{\gamma} \int_0^{\infty}\int_{\RR^{2d}} f^{\varepsilon}\partial_t\chi^{\varepsilon}\ud x\ud v\ud t -\varepsilon^{\gamma} \int_0^{\infty}\int_{\RR^{2d}} f_0(x,v)\chi^{\varepsilon}(0,x,v)\ud x\ud v \\
 &= \int_0^{\infty}\int_{\RR^{2d}} \big[ Q^+(f^{\varepsilon})\chi^{\varepsilon} - \nu f^{\varepsilon}\chi^{\varepsilon} + \varepsilon (v - j^\eps_F)  \cdot \nabla_x\chi^{\varepsilon}f^{\varepsilon}\big] \ud x\ud v\ud t  .
\end{align*}
Now, using \eqref{eq de chi^eps} and the fact that $Q^+(F) = \nu F$, we write
\begin{align*}
&- \int_0^{\infty}\int_{\RR^{2d}} f^{\varepsilon}\partial_t\chi^{\varepsilon}\ud x\ud v\ud t - \int_0^{\infty}\int_{\RR^{2d}} f_0(x,v)\chi^{\varepsilon}(0,x,v)\ud x\ud v \\
&= \varepsilon^{-\gamma}\int_0^{\infty}\int_{\RR^{2d}} \left(\big[f^{\varepsilon} (Q^+)^*(\chi^{\varepsilon}) - f^{\varepsilon}\nu \varphi \big] -\eps j^\eps_F  \cdot \nabla_x\chi^{\varepsilon}f^{\varepsilon} \right) \ud x\ud v\ud t \\
&= \varepsilon^{-\gamma}\int_0^{\infty}\int_{\RR^{2d}} \left( f^{\varepsilon} \big[(Q^+)^*(\chi^{\varepsilon}) - (Q^+)^*(\varphi) \big] -\eps j^\eps_F  \cdot \nabla_x\chi^{\varepsilon}f^{\varepsilon} \right) \ud x\ud v\ud t \\
&= \varepsilon^{-\gamma}\int_0^{\infty}\int_{\RR^{2d}} \big[Q^+(f^{\varepsilon}) \big(\chi^{\varepsilon} - \varphi \big) -\eps j^\eps_F  \cdot \nabla_x\chi^{\varepsilon}f^{\varepsilon} \big] \ud x\ud v\ud t ,
\end{align*}
and using the decomposition $f^\eps = \rho^\eps F + g^\eps$ with $g^\eps := f^\eps - \rho^\eps F$, we get
\begin{align*}
&- \int_0^{\infty}\int_{\RR^{2d}} f^{\varepsilon}\partial_t\chi^{\varepsilon}\ud x\ud v\ud t - \int_0^{\infty}\int_{\RR^{2d}} f_0(x,v)\chi^{\varepsilon}(0,x,v)\ud x\ud v   \\
&= \varepsilon^{-\gamma}\int_0^{\infty}\int_{\RR^{2d}} \left[ Q^+(\rho^{\varepsilon} F) \big(\chi^{\varepsilon} - \varphi \big) -\eps (j^\eps_F  \cdot \nabla_x \varphi) \rho^\eps F \right]\ud x\ud v\ud t \\
& \quad + \varepsilon^{-\gamma}\int_0^{\infty}\int_{\RR^{2d}} \left[Q^+(g^{\varepsilon}) \big(\chi^{\varepsilon} - \varphi \big) -\eps j^\eps_F  \cdot \nabla_x(\chi^\eps-\varphi)\rho^\eps F -\eps j^\eps_F  \cdot \nabla_x\chi^\eps g^\eps \right] \ud x\ud v\ud t .
\end{align*}
The second step of the proof of the main theorem consists in passing to the limit in the previous weak formulation. First, we rewrite it as follows:
\begin{align}
&- \int_0^{\infty}\int_{\RR^{2d}} f^{\varepsilon}\partial_t\chi^{\varepsilon}\ud x\ud v\ud t - \int_0^{\infty}\int_{\RR^{2d}} f_0(x,v)\chi^{\varepsilon}(0,x,v)\ud x\ud v  \label{donnees initiales} \\
&=   \varepsilon^{-\gamma}\int_0^{\infty}\int_{\RR^{2d}} \left[Q^+(g^{\varepsilon}) \big(\chi^{\varepsilon} - \varphi \big) -\eps j^\eps_F  \cdot \nabla_x(\chi^\eps-\varphi)\rho^\eps F -\eps j^\eps_F  \cdot \nabla_x\chi^\eps g^\eps \right] \ud x\ud v\ud t  \label{partie negligeable2}  \\
& \quad +  \varepsilon^{-\gamma}\int_0^{\infty}\int_{\RR^{d}}\rho^{\varepsilon} \int_{\RR^{d}}\nu F \left(\chi^{\varepsilon} - \varphi - \frac{\eps}{\nu} j^\eps_F \cdot \na_x \varphi \right) \ud v\ud x \ud t .  \label{partie importante}
\end{align}
The terms in \eqref{donnees initiales} and \eqref{partie negligeable2} are addressed in Section \ref{section 3.1}, while the terms in \eqref{partie importante} are dealt with separately in Section \ref{section 3.2}.

\subsection{First part of the limiting process}\label{section 3.1}

\subsubsection*{Limit of the time derivative}
Here, we deal with the limit of the integrals in \eqref{donnees initiales}.
\begin{lemma}
We have the following limits
$$ - \int_0^{\infty}\int_{\RR^{2d}} f^{\varepsilon}\partial_t\chi^{\varepsilon} \ \ud x\ud v\ud t \underset{\varepsilon \rightarrow 0}{\longrightarrow} - \int_0^{\infty}\int_{\RR^{d}} \rho\partial_t\varphi \ \ud x \ud t  $$
and 
$$- \int_{\RR^{d}} f_0\chi^{\varepsilon}(0,x,v) \ \ud x\ud v \underset{\varepsilon \rightarrow 0}{\longrightarrow} - \int_{\RR^{d}} \rho_0\varphi(0,x) \ \ud x  .$$
\end{lemma}

\ni The proof follows from Lemmas \ref{chi et phi} and \ref{lemme de compacite BL} and Lebesgue's theorem.  Indeed, we have
\begin{align*}
\bigg|\int_0^{\infty}\int_{\RR^{2d}} f^{\varepsilon}\big(\partial_t\chi^{\varepsilon}  -  \partial_t\varphi \big) \ud x\ud v\ud t \bigg| &\leqslant \ \| f^{\varepsilon} \|_{L^{\infty}(0,\infty,L^2_{F^{-1}}(\RR^{2d}))} \| \partial_t\chi^{\varepsilon} - \partial_t\varphi \|_{L^2(0,\infty,L^2_F(\RR^{2d}))} \\
&\leqslant \ \| f_0 \|_{L^2_{F^{-1}}(\RR^{2d})} \| \partial_t\chi^{\varepsilon} - \partial_t\varphi \|_{L^2(0,\infty,L^2_F(\RR^{2d}))} .
\end{align*}
By \eqref{F(chi_t - phi_t)-->0},  $\  \| \partial_t\chi^{\varepsilon} - \partial_t\varphi \|_{L^2(0,\infty,L^2_F(\RR^{2d}))}  \underset{\eps \to 0}{ \longrightarrow } 0 $.
Hence,  the limits of the previous lemma hold. 

\subsubsection*{Limit of a corrector terms}
The goal of this subsection is to show that the limit of \eqref{partie negligeable2} is zero. Note that there are more terms to handle here compared to \cite{M}, due to the presence of the drift term.
\begin{proposition}\label{partie negligeable}
For any test function $\varphi \in \mathcal{D}([0,\infty)\times\RR^{d})$,  let $\chi^{\varepsilon}$ be the function defined by \eqref{formule chi^eps}. Then,
$$ \lim_{\varepsilon\rightarrow 0} \ \varepsilon^{-\gamma}\int_0^{\infty}\int_{\RR^{2d}} \left[Q^+(g^{\varepsilon}) \big(\chi^{\varepsilon} - \varphi \big) -\eps j^\eps_F  \cdot \nabla_x\chi^\eps g^\eps -\eps j^\eps_F  \cdot \nabla_x(\chi^\eps-\varphi)\rho^\eps F \right] \ud x\ud v\ud t = 0 .  $$
\end{proposition}
\ni \bpp \ref{partie negligeable}.  The proof is given in three steps, where in each step we deal with a term, showing that it tends to 0. \smallskip

\ni \textbf{Step 1}. The aim of this step is to prove that, there exists a constant $C>0$ such that 
\begin{equation}\label{term1} 
	 \eps^{-\gamma}\left| \int_0^{\infty}\int_{\RR^{2d}} Q^+(g^{\varepsilon}) \left[\chi^{\varepsilon}(t,x,v) - \varphi(t,x) \right]  \ud v \ud x \ud t \right| \leqslant C  \big( \varepsilon^{\frac{2-\gamma}{2}} + \varepsilon^{\frac{\gamma}{2}}\big) \| \varphi \|_{H^1((0,\infty)\times\RR^{d})} .
\end{equation}
First,  using Assumption (B3), for all $t$ and $x$, we write 
\begin{align*}
\big| Q^+(g^\eps) \big| &\leqslant \bigg(\int_{\RR^{d}}\sigma^2(x,v,v')\frac{F(v')}{\nu(v')}\ud v' \bigg)^{\frac{1}{2}} \bigg( \int_{\RR^{d}}\big|g^\eps(v')\big|^2\frac{\nu(v')}{F(v')} \ud v'\bigg)^\frac{1}{2} \\
&= \nu F \bigg(\int_{\RR^{d}}\frac{\sigma^2}{F^2}\frac{F'}{\nu'}\frac{1}{\nu^2}\ud v' \bigg)^{\frac{1}{2}} \| g^\eps\|_{L^2_{\nu F^{-1}}(\RR^d)}  \\
&\leqslant M \nu(x,v) F(v) \| g\|_{L^2_{\nu F^{-1}}(\RR^d)} ,
\end{align*} 
where we used the notations $F':=F(v')$ and $\nu':=\nu(x,v')$. Also, we have thanks to Assumption (B2),  $\nu(x,v) \leqslant \nu_2 \langle v \rangle^\beta \leqslant C \lceil v \rfloor^{\beta}$. Therefore,
\begin{align*}
\bigg| \int Q^+(g^{\varepsilon})\big(\chi^{\varepsilon} - \varphi \big) \ud v \ud x \ud t \bigg| &\leqslant M \int_0^{\infty}\int_{\RR^{2d}} \nu F \big|\chi^{\varepsilon} - \varphi \big| \| g^\eps\|_{L^2_{\nu F^{-1}}(\RR^d)} \ud v \ud x \ud t \\
&\leqslant C \int_0^{\infty}\int_{\RR^{d}}\| g^\eps \|_{L^2_{\nu F^{-1}}(\RR^d)} \bigg(\int_{\RR^{d}} \lceil v \rfloor^{\beta} F \big|\chi^{\varepsilon} - \varphi \big| \ud v \bigg)  \ud x \ud t   \\
&\leqslant C \| g^\eps\|_{L^2_{\nu F^{-1}}((0,\infty)\times\RR^{2d})} \bigg[\int_0^{\infty}\int_{\RR^{d}} \bigg[\int_{\RR^{d}}  \lceil v \rfloor^{\beta}  F \big|\chi^{\varepsilon} - \varphi \big| \ud v \bigg]^2 \ud x \ud t\bigg]^\frac{1}{2}. 
\end{align*}
The constant $C$ may change from one line to another, but it does not depend on any parameter. Now, we have from  \eqref{g^eps < theta(eps) f_0}, $\| g^\eps\|_{L^2_{\nu F^{-1}}((0,\infty)\times\RR^{2d})} \leqslant C \eps^\frac{\gamma}{2} $.  Then 
\begin{equation}\label{int Q^+ g^eps < J_1+J_2}
\bigg| \int_0^{\infty}\int_{\RR^{2d}} Q^+(g^{\varepsilon})\big(\chi^{\varepsilon} - \varphi \big) \ \ud x\ud v\ud t \bigg| \leqslant C \eps^\frac{\gamma}{2}\big(J^\eps_1 + J^\eps_2 \big)^\frac{1}{2} ,
\end{equation}
where $$ J^\eps_1 := \int_0^{\infty}\int_{\RR^{d}}\bigg(\int_{|v| \leqslant \varepsilon^{-\frac{1}{1-\beta}}} \lceil v \rfloor^{\beta}F \big|\chi^{\varepsilon} - \varphi \big| \ud v \bigg)^2 \ud x \ud t $$
and $$ J^\eps_2 := \int_0^{\infty}\int_{\RR^{d}}\bigg(\int_{|v| \geqslant \varepsilon^{-\frac{1}{1-\beta}}} \lceil v \rfloor^{\beta} F \big|\chi^{\varepsilon} - \varphi \big| \ud v \bigg)^2 \ud x \ud t  .$$
\textbf{Estimation of $J^\eps_1$.} By decomposing the set $\{|v| \leqslant \varepsilon^{-\frac{1}{1-\beta}}\}$ into two parts $\{|v| \leqslant 1\}$ and $\{1\leqslant |v| \leqslant \varepsilon^{-\frac{1}{1-\beta}}\}$, we write
$$
J^\eps_1 \leqslant 2  \int_0^{\infty}\int_{\RR^{d}}\bigg[\bigg(\int_{|v| \leqslant 1} \lceil v \rfloor^{\beta}  F \big|\chi^{\varepsilon} - \varphi \big| \ud v  \bigg)^2 + \bigg( \int_{1 \leqslant |v| \leqslant \varepsilon^{-\frac{1}{1-\beta}}} \lceil v \rfloor^{\beta}  F \big|\chi^{\varepsilon} - \varphi \big| \ud v \bigg)^2 \bigg]\ud x \ud t . $$
Now let's deal with the integral
$$ \int_0^{\infty}\int_{\RR^{d}}\bigg(\int_{\Omega} \lceil v \rfloor^{\beta} F \big|\chi^{\varepsilon} - \varphi \big| \ud v \bigg)^2 \ud x \ud t  ,  $$
with  $\Omega = \{|v| \leqslant 1 \}$ or $\Omega = \{1 \leqslant |v| \leqslant \varepsilon^{-\frac{1}{1-\beta}}\}$.   As in the proof of Lemma \ref{chi et phi}, using Taylor's formula,  Assumption (B2) and the Cauchy-Schwarz inequality, we write
\begin{align*}
 \big| \chi^{\varepsilon} - \varphi \big| &= \bigg| \int_0^\infty \lceil v \rfloor^{-\beta} \tilde \nu_\eps(x,v,z) \ e^{-\int_0^z \lceil v \rfloor^{-\beta} \tilde \nu_\eps(x,v,s) \ud s} \big[\tilde \varphi_\eps(t,x,v,z) - \varphi(t,x) \big] \ud z \bigg|  \\
 &\leqslant \int_0^{\infty}\nu_2 e^{-\nu_1 z}\varepsilon |v| \lceil v \rfloor^{-\beta} z \int_0^1 |\na_x \varphi(t,x+\eps v \lceil v \rfloor^{-\beta}zs)| \ \ud s \ud z \\
 &\leqslant \nu_2\varepsilon |v| \lceil v \rfloor^{-\beta} \bigg(\int_0^{\infty}z^2 e^{-\nu_1 z}\ud z\bigg)^\frac{1}{2} \bigg(\int_0^{\infty} e^{-\nu_1 z}\int_0^1 \big| \na_x\varphi(t,x+\varepsilon |v|^{-\beta}vzs)\big|^2 \ud s\ud z\bigg)^\frac{1}{2} .
 \end{align*}
Therefore,
\begin{align*}
& \int_0^{\infty}\int_{\RR^{d}}\bigg(\int_{\Omega} \lceil v \rfloor^{\beta} F \big|\chi^{\varepsilon} - \varphi \big| \ud v \bigg)^2 \ud x \ud t   \\ 
&\leqslant C \int_0^{\infty}\int_{\RR^{d}} \bigg[\int_{\Omega} \varepsilon |v| F\bigg(\int_0^{\infty} e^{-\nu_1 z}\int_0^1 \big| \na_x\varphi(t,x+\varepsilon |v|^{-\beta}vzs)\big|^2 \ud s\ud z\bigg)^\frac{1}{2}\ud v\bigg]^2\ud x \ud t  \\
&\leqslant C\varepsilon^2\bigg( \int_{\Omega} |v| F\ud v\bigg) \bigg(\int_0^{\infty}\int_{\RR^{d}}\int_{\Omega} |v| F \int_0^{\infty} e^{-\nu_1 z}\int_0^1 \big| \na_x\varphi(t,x+\varepsilon |v|^{-\beta}vzs)\big|^2 \ud s\ud z\ud v \bigg) \ud x \ud t ,
\end{align*}
where $C := \nu_2 \int_0^{\infty}z^2 e^{-\nu_1 z}\ud z $.  Now, taking $\Omega = \{|v| \leqslant 1 \}$ in the previous inequality, we obtain
\begin{align*}
&\int_0^{\infty}\int_{\RR^{d}}\bigg(\int_{|v| \leqslant 1} \lceil v \rfloor^{\beta}  F \big|\chi^{\varepsilon} - \varphi \big| \ud v \bigg)^2 \ud x \ud t \\
&\leqslant C\varepsilon^2 \bigg( \int_{|v| \leqslant 1} |v| F\ud v\bigg)^2  \int_0^{\infty}e^{-\nu_1 z} \int_0^1 \bigg(\int_0^{\infty}\int_{\RR^{d}}\big| \na_x\varphi(t,x+\varepsilon |v|^{-\beta}vzs)\big|^2\ud x \ud t \bigg)\ud s \ud z \\
&\leqslant C \varepsilon^2\| \na_x\varphi \|^2_{L^2((0,\infty)\times\RR^{d})} .
\end{align*} 
Similarly,  for $\Omega = \{1 \leqslant |v| \leqslant \varepsilon^{-\frac{1}{1-\beta}}\}$:
\begin{align*}
&\int_0^{\infty}\int_{\RR^{d}}\bigg(\int_{1 \leqslant |v| \leqslant \varepsilon^{-\frac{1}{1-\beta}}} |v|^{\beta}  F \big|\chi^{\varepsilon} - \varphi \big| \ud v \bigg)^2 \ud x \ud t \\ 
&\leqslant C\varepsilon^2 \bigg( \int_{1 \leqslant |v| \leqslant \varepsilon^{-\frac{1}{1-\beta}}} |v| F\ud v\bigg)^2  \int_0^{\infty}e^{-\nu_1 z} \int_0^1 \bigg(\int_0^{\infty}\int_{\RR^{d}} \big| \na_x\varphi(t,x+\varepsilon |v|^{-\beta}vzs)\big|^2\ud x \ud t \bigg)\ud s \ud z \\
& \leqslant  C \varepsilon^2 \| \na_x\varphi \|^2_{L^2((0,\infty)\times\RR^{d})} \bigg(\int_{1 \leqslant |v| \leqslant \varepsilon^{-\frac{1}{1-\beta}}} |v| F \ud v\bigg)^2  \\
&\leqslant C\varepsilon^2 \| \na_x\varphi \|^2_{L^2((0,\infty)\times\RR^{d})} \big(1+\eps^{\gamma-1}\big)^2 ,
\end{align*}
where we used the inequality $F(v) \leqslant C |v|^{-\alpha-d}$ for $|v| \geqslant 1$, thanks to Assumption (B1), and where we recall that $\gamma:=\frac{\alpha-\beta}{1-\beta}$.
Hence,
\begin{equation}\label{estimation J^eps_1 BL}
J^\eps_1 \leqslant C \| \na_x\varphi \|^2_{L^2((0,\infty)\times\RR^{d})}\big(\varepsilon^\gamma +\varepsilon \big)^2  .
\end{equation}
\textbf{Estimation of $J^\eps_2$.} For the range of velocities $\{|v| \geqslant \varepsilon^{-\frac{1}{1-\beta}}\}$,  we write
\begin{align*}
J^\eps_2 &:= \int_0^{\infty}\int_{\RR^{d}}\bigg(\int_{|v| \geqslant \varepsilon^{-\frac{1}{1-\beta}}}  |v|^{\beta}  F \big|\chi^{\varepsilon} - \varphi \big| \ud v \bigg)^2\ud x \ud t \\
&\leqslant \bigg(\int_{|v| \geqslant \varepsilon^{-\frac{1}{1-\beta}}} |v|^{\beta} F\ud v\bigg)\bigg(\int_0^{\infty}\int_{\RR^{d}}\int_{|v| \geqslant \varepsilon^{-\frac{1}{1-\beta}}} |v|^{\beta} F \big|\chi^{\varepsilon} - \varphi \big|^2 \ud v \ud x\ud t \bigg) \\
&\leqslant 2 \bigg(\int_{|v| \geqslant \varepsilon^{-\frac{1}{1-\beta}}} |v|^{\beta} F\ud v\bigg) \int_{|v| \geqslant \varepsilon^{-\frac{1}{1-\beta}}} |v|^{\beta} F \bigg( \int_0^{\infty}\int_{\RR^{d}} \big(  |\chi^\varepsilon|^2 + | \varphi |^2\big)\ud x \ud t \bigg)\ud v .
\end{align*}
Since $\| \chi^\varepsilon \|_{L^2_{x,t}} \leqslant \frac{\nu_2}{\nu_1}\| \varphi \|_{L^2}$ (see proof of inequality \eqref{chi_t < phi_t}),  and since we have $F(v) \leqslant C |v|^{-\alpha-d}$ for $|v| \geqslant \varepsilon^{-\frac{1}{1-\beta}}$ then,
\begin{equation}\label{estimation J^eps_2 BL}
J^\eps_2 \leqslant C \| \varphi \|^2_{L^2((0,\infty)\times\RR^{d})}\bigg(\int_{|v| \geqslant \varepsilon^{-\frac{1}{1-\beta}}} |v|^{\beta} F \ud v\bigg)^2  \leqslant C  \| \varphi \|^2_{L^2((0,\infty)\times\RR^{d})} \varepsilon^{2 \gamma}.  
\end{equation}
Finally, going back to \eqref{int Q^+ g^eps < J_1+J_2} and using estimates \eqref{estimation J^eps_1 BL} and \eqref{estimation J^eps_2 BL}, we obtain
$$ \bigg| \int_0^{\infty}\int_{\RR^{2d}} Q^+(g^{\varepsilon})\big(\chi^{\varepsilon} - \varphi \big) \ud v \ud x \ud t  \bigg| \leqslant C\varepsilon^{\frac{\gamma}{2}}\big( \varepsilon + \varepsilon^{\gamma}\big) \| \varphi \|_{H^1((0,\infty)\times\RR^{d})}.$$
Which implies that,
$$   \bigg| \varepsilon^{-\gamma} \int_0^{\infty}\int_{\RR^{2d}} Q^+(g^{\varepsilon})\big(\chi^{\varepsilon} - \varphi \big) \ud v \ud x \ud t  \bigg| \leqslant C \big( \varepsilon^{\frac{2-\gamma}{2}} + \varepsilon^{\frac{\gamma}{2}}\big)\| \varphi \|_{H^1((0,\infty)\times\RR^{d})} .$$
Hence inequality \eqref{term1} holds. Moreover, since $\gamma \in(0,2)$ we get
 $$ \underset{\eps \to 0}{\lim} \ \eps^{-\gamma} \int_0^{\infty}\int_{\RR^{2d}} Q^+(g^{\varepsilon})\big[\chi^{\varepsilon}(t,x,v) - \varphi(t,x) \big] \ud v \ud x \ud t = 0 .$$
\textbf{Step 2.} The aim of this step is to prove that, there exists a constant $C>0$ such that
$$ \eps^{1-\gamma} \left| \int_0^{\infty}\int_{\RR^{2d}} j^\eps_F  \cdot \nabla_x \chi^\eps g^\eps  \ud v \ud x \ud t \right| \leqslant C \eps^{\min\{\frac{2-\gamma}{2};\frac{1}{1-\beta}\}} \|j^\eps_F\|_{L^\infty} \|\varphi\|_{H^1} .  $$
First, to avoid long expressions, we introduce the notation $\bar{\nu}_\eps := \lceil v \rfloor^{-\beta} \tilde{\nu}(x+\eps \lceil v \rfloor^{-\beta} vz,v)$, and we recall that $\bar{\nu}_\eps$ is bounded below and above by $\nu_1$ and $\nu_2$ thanks to (B2). For simplicity, we shall not keep track of all constants and simply write that $\mathrm A \lesssim \mathrm B$ if there is a positive constant $C>0$ such that, $\mathrm A \leqslant C \mathrm B$. Now, we have
$$ \na_x \chi^\eps = \int_0^\infty \left[ \bigg(\na_x\bar{\nu}_\eps -\bar{\nu}_\eps \int_0^z \bar{\nu}_\eps \ud s \bigg) \tilde \varphi + \bar{\nu}_\eps \na_x \tilde \varphi_\eps \right] e^{-\int_0^z \bar{\nu}_\eps \ud s}   \ud z .$$
Then, using the two inequalities from (B2) on $\nu$ and $\na_x \nu$, we obtain 
$$ \int \left| j^\eps_F \cdot \na_x \chi^\eps g^\eps \right| \ud v \ud x \ud t \lesssim \|j^\eps_F\|_{L^\infty} \int \int_0^\infty (1+z) \big(|\tilde \varphi_\eps| + |\na_x \tilde \varphi_\eps|\big) e^{-\nu_2 z} \ud z \ |g^\eps| \ \ud x\ud v\ud t .$$
Now, we distinguish two cases, $\beta \geqslant 0$ and $\beta <0$. For $\beta \geqslant 0$, by Cauchy-Schwarz inequqlity 
\begin{align*}
	\int &\left| j^\eps_F  \cdot \na_x \chi^\eps g^\eps \right| \ud v \ud x \ud t \\
	&\lesssim \|j^\eps_F\|_{L^\infty} \|g^\eps\|_{L^2_{\nu F^{-1}}} \bigg[\int \bigg(\int_0^\infty (1+z) \big(|\tilde \varphi_\eps| + |\na_x \tilde \varphi_\eps|\big) e^{-\nu_2 z} \ud z\bigg)^2 \frac{F}{\nu}  \ud x\ud v\ud t \bigg]^{\frac{1}{2}} \\
	&\lesssim \|j^\eps_F\|_{L^\infty} \|g^\eps\|_{L^2_{\nu F^{-1}}} \| \varphi \|_{H^1} \bigg(\int_0^\infty (1+z)^2  e^{-\nu_2 z} \ud z  \int \frac{F}{\nu} \ud v \bigg)^{\frac{1}{2}}  \\
	&\lesssim \eps^{\frac{\gamma}{2}} \|j^\eps_F\|_{L^\infty((0,\infty)\times\RR^{2d})} \| \varphi \|_{H^1((0,\infty)\times\RR^{d})} ,
\end{align*}
since $\|g^\eps\|_{L^2_{\nu F^{-1}}} \lesssim \eps^{\frac{\gamma}{2}}$ and $\int \frac{F}{\nu} \ud v \lesssim \int \langle v \rangle^{-\alpha-\beta-d} \ud v \lesssim 1$ for $\alpha >0$ and $\beta \geqslant 0$. \smallskip

\ni For $\beta <0$, we need to split the integral over $v$ into two parts, $\Omega_\eps := \{ |v| \leqslant \eps^{-\frac{1}{1-\beta}} \}$ and $\Omega_\eps^c$. Thus, by Cauchy-Schwarz inequqlity we write
\begin{align*}
	&\int \left| j^\eps_F \cdot \na_x \chi^\eps g^\eps \right| \ud v \ud x \ud t \\
	&\lesssim  \|j^\eps_F\|_{L^\infty} \bigg[ \bigg(\int |g^\eps|^2\frac{\nu}{F}\ud v \ud x \ud t\bigg)^{\frac{1}{2}}\bigg(\int_{t,x} \int_{\Omega_\eps} \bigg[\int_0^\infty (1+z) \big(|\tilde \varphi_\eps| + |\na_x \tilde \varphi_\eps|\big) e^{-\nu_2 z} \ud z\bigg]^2 \frac{F}{\nu} \ud v \ud x\ud t\bigg)^{\frac{1}{2}} \\
	& \quad +\bigg(\int_{t,x} \int_{\Omega_\eps
	^c} |g^\eps|^2\frac{1}{F}\ud v \ud x \ud t\bigg)^{\frac{1}{2}}\bigg(\int_{t,x} \int_{\Omega_\eps^c} \bigg[\int_0^\infty (1+z) \big(|\tilde \varphi_\eps| + |\na_x \tilde \varphi_\eps|\big) e^{-\nu_2 z} \ud z\bigg]^2 \ F \ \ud v \ud x\ud t\bigg)^{\frac{1}{2}} \bigg].
\end{align*}
Now, since $\nu$ behaves like $|v|^\beta$ and $\|g^\eps\|_{L^2_{\nu F^{-1}}} \lesssim \eps^{\frac{\gamma}{2}}$ with $\gamma := \frac{\alpha-\beta}{1-\beta}$, then 
$$ \int_{t,x} \int_{\Omega_\eps^c} |g^\eps|^2\frac{1}{F}\ud v \ud x \ud t \lesssim \eps^\frac{\beta}{1-\beta} \int_0^{\infty}\int_{\RR^{2d}} |g^\eps|^2\frac{\nu}{F}\ud v \ud x \ud t \lesssim \eps^{\frac{\alpha}{1-\beta}} .$$
	Thus, the previous inequality, together with 
	$$\int_{\Omega_\eps} \frac{F}{\nu} \ \ud v \lesssim 1+\eps^{\frac{\alpha+\beta}{1-\beta}} \quad \mbox{ and } \quad \int_{\Omega_\eps^c} F \ \ud v \lesssim \eps^{\frac{\alpha}{1-\beta}}$$
	imply that
	$$ \int \left| j^\eps_F \cdot \na_x \chi^\eps g^\eps \right| \ud v \ud x \ud t \lesssim \eps^{\frac{1}{1-\beta}} \|j^\eps_F\|_{L^\infty((0,\infty)\times\RR^{2d})} \| \varphi \|_{H^1((0,\infty)\times\RR^{d})}  .$$
	Finally, for all $\beta < 1$
	$$ \eps^{-\gamma} \left| \int \eps j^\eps_F \cdot \na_x \chi^\eps g^\eps \ \ud v \ud x \ud t \right|\lesssim \eps^{\min\{1-\frac{\gamma}{2};\frac{1}{1-\beta}\}} \|j^\eps_F\|_{L^\infty}  \|\varphi\|_{H^1}  . $$	
Consequently,
$$ \lim_{\varepsilon\rightarrow 0} \ \eps^{-\gamma}\int_0^{\infty}\int_{\RR^{2d}} \eps j^\eps_F  \cdot \nabla_x\chi^\eps g^\eps \ \ud v \ud x \ud t = 0 .  $$
\textbf{Step 3.} The aim of this step is to prove that, there exists a constant $C>0$ such that
$$ \eps^{1-\gamma} \left|\int_0^{\infty}\int_{\RR^{2d}} j^\eps_F  \cdot \nabla_x(\chi^\eps-\varphi)\rho^\eps F \ \ud v \ud x \ud t \right| \leqslant C  \big(\eps^{2-\gamma} +\eps^\frac{1}{1-\beta} \big) \|j^\eps_F\|_{L^\infty} \|\varphi\|_{L^1_t H^2_x} \|\rho^\eps\|_{L^\infty_t L^2_x}. $$
The proof of this inequality is given in the same spirit as the previous step. Indeed, by the Cauchy-Schwarz inequality with respect to $x$, and since $\rho^\eps$ depends only on $t$ and $x$, and is uniformly bounded in $L^{\infty}_t L^2_x$ by Lemma \ref{lemme de compacite BL}, we have
\begin{align*}
\left|\int j^\eps_F  \cdot \nabla_x(\chi^\eps-\varphi)\rho^\eps F \ \ud v \ud x \ud t \right| &\leqslant \|j^\eps_F\|_{L^\infty} \|\rho^\eps\|_{L^{\infty}_t L^2_x} \int \bigg[\int \bigg(\int \big|\na_x(\chi^\eps-\varphi)\big|F \ud v\bigg)^2 \ud x \bigg]^\frac{1}{2}\ud t .
\end{align*}
It remains to estimate the integral on the right-hand side. With the notation of the previous step, one has
\begin{align}\label{grad chi-phi}
\left|\na_x (\chi^\eps-\varphi) \right| &= \left| \int_0^\infty \left[ \bigg(\na_x\bar{\nu}_\eps -\bar{\nu}_\eps \int_0^z \bar{\nu}_\eps \ud s \bigg) (\tilde \varphi_\eps-\varphi) + \bar{\nu}_\eps \na_x (\tilde \varphi_\eps-\varphi) \right] e^{-\int_0^z \bar{\nu}_\eps \ud s}   \ud z \right| \nonumber \\
&\lesssim \int_0^\infty (1+z) \big(|\tilde \varphi_\eps-\varphi| + |\na_x (\tilde \varphi_\eps-\varphi)| \big) e^{-\nu_2 z} \ud z .
\end{align}
Now, we decompose the integral over $v$ into two parts, $\Omega_\eps$ and $\Omega_\eps^c$, as follows:
$$ \int_{\RR^d} \big|\na_x(\chi^\eps-\varphi)\big|F \ud v = \int_{\Omega_\eps} \big|\na_x(\chi^\eps-\varphi)\big|F \ud v + \int_{\Omega_\eps^c} \big|\na_x(\chi^\eps-\varphi)\big|F \ud v . $$
On $\Omega_\eps$, using a Taylor expansion, we write
\begin{align*} 
|\tilde \varphi_\eps-\varphi| &= \left| \int_0^1 \eps \lceil v \rfloor^{-\beta} z v \cdot \na_x\varphi(t,x+\varepsilon \lceil v \rfloor^{-\beta}vzs)\ud s  \right| \\
&\lesssim \eps \langle v \rangle^{1-\beta} z \bigg( \int_0^1 \big|\na_x\varphi(t,x+\eps \lceil v \rfloor^{-\beta}vzs)\big|^2 \ud s\bigg)^\frac{1}{2} .
\end{align*}
Similarly for the gradient. This implies that 	
\begin{align*}
&\left|\na_x (\chi^\eps-\varphi) \right|
\\  & \qquad \lesssim \eps \langle v \rangle^{1-\beta} \int_0^\infty (z+z^2) e^{-\nu_2 z}\bigg(\int_0^1 \big(\big|\na_x\tilde \varphi_\eps(t,x,v,zs)\big|^2 + \big|\na_x^2\tilde \varphi_\eps(t,x,v,zs)\big|^2\big) \ud s\bigg)^\frac{1}{2}  \ud z \\
& \qquad \lesssim \eps \langle v \rangle^{1-\beta} \bigg(\int_0^\infty \int_0^1 \big(\big|\na_x\tilde \varphi_\eps(t,x,v,zs)\big|^2 + \big|\na_x^2\tilde \varphi_\eps(t,x,v,zs)\big|^2\big)\ud s \ e^{-\nu_2 z} \ud z \bigg)^\frac{1}{2}.
\end{align*}
Therefore,
$$ \int \bigg(\int_{\Omega_\eps} \big|\na_x(\chi^\eps-\varphi)\big|F \ud v\bigg)^2 \ud x \lesssim \eps^2 \bigg(\int_{\Omega_\eps} \langle v \rangle^{1-\beta} F \ud v \bigg)^2 \big(\|\na_x \varphi(t)\|_{L^2_x}^2 + \|\na_x^2 \varphi(t)\|_{L^2_x}^2\big). $$
Now, since $\int_{\Omega_\eps} \langle v \rangle^{1-\beta} F \ud v \lesssim 1+\eps^{\frac{\alpha+\beta-1}{1-\beta}}$ then,
$$ \int_0^\infty \bigg[\int_{\RR^d} \bigg(\int_{\Omega_\eps} \big|\na_x(\chi^\eps-\varphi)\big|F \ud v\bigg)^2 \ud x \bigg]^\frac{1}{2}\ud t \lesssim \big(\eps+\eps^{\frac{\alpha}{1-\beta}}\big) \ \|\varphi\|_{L^1_t H^2_x}. $$
On $\Omega_\eps^c$, using inequality \eqref{grad chi-phi}, applying the Cauchy-Schwarz inequality in $z$, and integrating over $v$ and $x$, we obtain 
$$ \int \bigg(\int_{\Omega_\eps^c} \big|\na_x(\chi^\eps-\varphi)\big|F \ud v \bigg)^2 \ud x \lesssim \bigg(\int_{\Omega_\eps^c} F \ud v \bigg)^2 \big(\|\varphi(t)\|_{L^2_x}^2 + \|\na_x \varphi(t)\|_{L^2_x}^2\big).$$
From where,
$$ \int_0^\infty \bigg[\int_{\RR^d} \bigg(\int_{\Omega_\eps^c} \big|\na_x(\chi^\eps-\varphi)\big|F \ud v\bigg)^2 \ud x \bigg]^\frac{1}{2}\ud t \lesssim \eps^{\frac{\alpha}{1-\beta}} \ \|\varphi\|_{L^1_t H^1_x}. $$
Hence the proof of the inequality of this step is complete. As a consequence, we have the following limit
$$ \lim_{\varepsilon\rightarrow 0} \ \eps^{-\gamma}\int_0^{\infty}\int_{\RR^{2d}} \eps j^\eps_F  \cdot \nabla_x(\chi^\eps-\varphi)\rho^\eps F \ \ud v \ud x \ud t = 0 .  $$
Thus, by combining the limits obtained in the three steps, we get the limit of Proposition \ref{partie negligeable}, and the proof of this proposition is complete.
\epp

To conclude the proof of the main theorem, it remains to pass to the limit in the last line of the weak formulation, integrals \eqref{partie importante}, which is the subject of the following subsection.
\subsection{Obtention of the limiting operator}\label{section 3.2}
For the passage to the limit in \eqref{partie importante}, we have the following proposition.
\begin{proposition}\label{convergence vers Laplacien}
For any test function $\varphi \in \mathcal{D}([0,\infty)\times\RR^{d})$,  let $\chi^{\varepsilon}$ the function defined by \eqref{formule chi^eps}. Then,
$$ \lim_{\eps\to 0} \ \eps^{-\gamma}\int_{\RR^{d}}\nu(x,v)F(v)\left[\chi^{\eps}(t,x,v) - \varphi(t,x) -\frac{\eps}{\nu(x,v)} j^\eps_F \cdot \na_x \varphi(t,x)\right] \ud v = -\kappa \mathcal{L}^*(\varphi) ,$$
where $\mathcal{L}$ is the operator defined in Theorem \ref{main BL}, formula \eqref{definition de L BL}.  Moreover, the last limit is uniform with respect to $x$ and $t$.
\end{proposition}

\ni \bpp \ref{convergence vers Laplacien}.  The proof is done in two main steps. The first consists in showing that the small velocities do not participate in the limit, and therefore it will be the large velocities which gives the elliptic operator $\mathcal{L}$,  which is the subject of the second step.    
 
\ni Before doing this, we need to distinguish three cases according to the different expressions of the drift $j^\eps_F$, and rewrite the integral in the limit of the proposition above in a more compact form. Indeed, since $j^\eps_F$ and $\varphi$ do not depend on $v$ and $\int_{\mathbb{R}^d} F \ \ud v = 1$, and thanks to the identity
$$ \int_0^\infty z \lceil v \rfloor^{-\beta} \nu (x,v) e^{- \lceil v \rfloor^{-\beta}  \nu(x,v) z}  \ud z = \frac{1}{\lceil v \rfloor^{-\beta} \nu (x,v)}  ,$$
the integral corresponding to the drift term can be written as follows: \\
$\bullet$ For $\alpha > 1$, we have 
$$ \int_{\RR^{d}} \nu F \left(\frac{\eps}{\nu} j^\eps_F \cdot \na_x \varphi \right) \ud v = \eps j^\eps_F \cdot \na_x \varphi = \int_{\RR^{d}} \nu F \int_0^\infty \left(\eps \lceil v \rfloor^{-\beta}  v z\cdot \na_x \varphi \right) \lceil v \rfloor^{-\beta}  \nu \ e^{- \lceil v \rfloor^{-\beta}  \nu z} \ud z \ud v .$$
$\bullet$ For $\alpha = 1$, we have with $\Omega_\eps := \{|v| \leqslant \eps^{-\frac{1}{1-\beta}}\}$
$$ \int_{\RR^{d}} \nu F \left(\frac{\eps}{\nu} j^\eps_F \cdot \na_x \varphi \right) \ud v = \eps j^\eps_F \cdot \na_x \varphi = \int_{\Omega_\eps} \nu F \int_0^\infty \left(\eps \lceil v \rfloor^{-\beta}  v z\cdot \na_x \varphi \right) \lceil v \rfloor^{-\beta}  \nu \ e^{- \lceil v \rfloor^{-\beta}  \nu z}  \ud z \ud v .$$
$\bullet$ For $\alpha < 1$, the integral vanishes since $j^\eps_F=0$. \smallskip

\ni Finally, we intoduce two functions $P_\eps$ and $P_0$ defined  by
$$
P_{\varepsilon}(x,v,z) := \lceil v \rfloor^{-\beta} \tilde \nu_\eps(x,v,z)\ e^{-\int_0^z \lceil v \rfloor^{-\beta} \tilde \nu_\eps(x,v,s) \ud s} \ \mbox{and} \  P_0(x,v,z) := \lceil v \rfloor^{-\beta} \nu(x,v) \ e^{ -z \lceil v \rfloor^{-\beta} \nu(x,v)} .  $$
It is in this part that the drift term $j^\eps_F$ plays a very important role, since the symmetry of $F$ and $\nu$ is not assumed. In comparison to the proof given in \cite{M}, symmetry allows one to recover the identity  
$$ \nu(x,-v)F(-v)P_0(x,-v,z) = \nu(x,v)F(v)P_0(x,v,z) , $$ 
which is not guaranteed for us. This latter identity allows, through a second-order Taylor expansion, to gain a factor of $\eps$ for bounded $|v|$. With the presence of the drift, we rewrite the integral in Proposition \ref{convergence vers Laplacien} as follows: \\
$\bullet$ For $\alpha > 1$, we have
\begin{equation*}
\int_{\RR^{d}} \nu F\left(\chi^{\varepsilon}-\varphi - \frac{\eps}{\nu} \ j^\eps_F \cdot \na_x\varphi \right) \ud v  = \int_{\RR^{d}} \nu F \int_0^\infty \left( P_{\varepsilon}  [\tilde \varphi_\eps -\varphi] - P_0 \eps \lceil v \rfloor^{-\beta}  v z\cdot \na_x \varphi \right) \ud z \ud v .
\end{equation*}
$\bullet$ For $\alpha = 1$, we have
\begin{align*}
\int_{\RR^{d}} \nu F\left(\chi^{\varepsilon}-\varphi - \frac{\eps}{\nu} \ j^\eps_F \cdot \na_x\varphi \right) \ud v  =& \int_{\Omega_\eps} \nu F \int_0^\infty \left( P_{\varepsilon}  [\tilde \varphi_\eps -\varphi] - P_0 \eps \lceil v \rfloor^{-\beta}  v z\cdot \na_x \varphi \right) \ud z \ud v \\
&+\int_{\Omega_\eps^c} \nu F \int_0^\infty P_{\varepsilon}[\tilde \varphi_\eps -\varphi] \ud z \ud v. 
\end{align*}
$\bullet$ For $\alpha < 1$, we have
\begin{equation}\label{petites vitesses alpha<1}
\int_{\RR^{d}} \nu F\left(\chi^{\varepsilon}-\varphi - \frac{\eps}{\nu} \ j^\eps_F \cdot \na_x\varphi \right) \ud v  = \int_{\RR^{d}} \nu F \int_0^\infty P_{\varepsilon} [\tilde \varphi_\eps -\varphi] \ud z \ud v .
\end{equation}

\begin{remark}
	In order to simplify the computations and without loss of generality, we can replace Assumptions (B1) and (B2) concerning the behavior of $F$ and $\nu$ at infinity by
$$ F(v) = \frac{\kappa}{|v|^{\alpha+d}} \quad \mbox{ and } \quad \nu(x,v) = \nu_0(x) |v|^{-\beta} ,\quad \mbox{for all } \ |v| \geqslant  C  \ \mbox{ and } \ x \in \RR^d.$$
\end{remark}

\ni \textbf{Step 1: Small velocities don't contribute to the limit}
\begin{lemma}\label{petites vitesses BL}
Let $ \varphi \in \mathcal{D}([0,+\infty)\times\RR^d)$ and $\chi^\eps$ defined by \eqref{formule chi^eps}.  Assume (A1-A2) and (B2).  Then, the following uniforme estimate holds
$$\bigg| \int_{|v|\leqslant C} \nu F \left[\chi^{\varepsilon}-\varphi -\eps \frac{v}{\nu}\cdot\na_x \varphi \right] \ud v \bigg| \leqslant C_0 \eps^2 \| \varphi \|_{W^{2,\infty}} ,$$
where $C_0 > 0$ depends on $\nu_1$,  $\nu_2$ and $C$, the constant of Assumption (B2).  Consequently,
$$ \lim_{\varepsilon\rightarrow 0} \ \varepsilon^{-\gamma}\int_{|v|\leqslant C}\nu(x,v)F(v)\left[\chi^{\varepsilon}(t,x,v) - \varphi(t,x) - \frac{\eps}{\nu(x,v)} \ j^\eps_F \cdot \na_x \varphi(t,x)\right] \ud v = 0 .$$
\end{lemma}
\ni \bpl \ref{petites vitesses BL}. For $\alpha \geqslant 1$, we have
\begin{align*}
\int_{|v|\leqslant C} \nu F \left[\chi^{\varepsilon}-\varphi -\eps \frac{v}{\nu}\cdot\na_x \varphi \right] \ud v =& \int_{|v|\leqslant C} \nu F\int_0^\infty \big( P_\varepsilon - P_0 \big) \big[\tilde \varphi_\eps - \varphi \big] \ud z \ud v \\
&+ \int_{|v|\leqslant C} \nu F \int_0^\infty  P_0 \big[\tilde \varphi_\eps - \varphi - \varepsilon v \lceil v \rfloor^{-\beta}z\cdot \nabla_x\varphi \big] \ud z \ud v .
\end{align*}
Let's start with the integral of the first line.  We have
\begin{align*}
\big| P_\varepsilon - P_0 \big| =& \bigg|  \lceil v \rfloor^{-\beta} \tilde \nu_\eps(x,v,z)\ e^{ -\int_0^z \lceil v \rfloor^{-\beta} \tilde \nu_\eps(x,v,s) \ud s} -  \lceil v \rfloor^{-\beta} \nu(x,v)\ e^{-z  \lceil v \rfloor^{-\beta}\nu(x,v)} \bigg| \\
\leqslant &  \lceil v \rfloor^{-\beta} \big| \tilde \nu_\eps(x,v,z) - \nu(x,v) \big|e^{ -\int_0^z  \lceil v \rfloor^{-\beta} \tilde \nu_\eps(x,v,s) \ud s} \\
& \ + \lceil v \rfloor^{-\beta} \nu(x,v)\bigg|e^{-\int_0^z  \lceil v \rfloor^{-\beta} |\tilde \nu_\eps(x,v,s)-\nu(x,v)| \ud s} -1\bigg| e^{ -\int_0^z \lceil v \rfloor^{-\beta} \min\{\tilde\nu_\eps(x,v,s),\nu(x,v)\}\ud s}  .
\end{align*}
Now,  since $\ |e^{-y} - 1| \leqslant y \ $  for all $ \ y \geqslant 0$, then
$$  \bigg|e^{-\int_0^z  \lceil v \rfloor^{-\beta} |\tilde \nu_\eps(x,v,s)-\nu(x,v)| \ud s} -1\bigg|  \leqslant \int_0^z \lceil v \rfloor^{-\beta} \big|\tilde \nu_\eps(x,v,s)-\nu(x,v) \big| \ud s.$$
And since $\nu$ is $C^1$ with respect to $x$ and $\| \langle v\rangle^{-\beta} \pa_x \nu \|_{L^\infty} \leqslant C$, then
$$ \big|\tilde \nu_\eps(x,v,z)-\nu(x,v) \big| := \big|\nu(x,x+\varepsilon \lceil v \rfloor^{-\beta} vz,v)-\nu(x,v) \big| \leqslant \eps z |v| \| \lceil v \rfloor^{-\beta} \pa_x \nu \|_{L^\infty}  \leqslant C \eps z |v| . $$
Hence,
\begin{align*}
\big| P_\varepsilon - P_0 \big| &\leqslant  C \eps |v|  \lceil v \rfloor^{-\beta} z^2 e^{-\nu_1 z} .
\end{align*}
We have also,
\begin{equation}\label{estimation phi_eps -phi}
	\big| \tilde \varphi_\eps(t,x,v,z) - \varphi(t,x) \big|  := \big| \varphi(t,x+\varepsilon |v|^{-\beta}vz) - \varphi(t,x) \big| \leqslant C\eps z |v| \lceil v \rfloor^{-\beta} \| \na_x\varphi \|_{L^{\infty}}  .
\end{equation}
Thus,  on the one hand,
\begin{align*}
&\bigg| \int_{|v|\leqslant C} \nu F \int_0^\infty \big( P_\varepsilon - P_0 \big) \big[ \tilde \varphi_\eps(t,x,v,z) - \varphi(t,x) \big] \ud z \ud v \bigg| \\
&\leqslant C \int_{|v| \leqslant C} \lceil v \rfloor^{\beta}  F \int_0^\infty\varepsilon ^2 e^{-\nu_1 z} |v|^{2} \lceil v \rfloor^{-2\beta}(z+z^2)\| \na_x\varphi\|_{L^{\infty}} \ud z\ud v \\
&\leqslant C\varepsilon^2 \| \na_x\varphi\|_{L^{\infty}} .
\end{align*}
On the other hand,
\begin{align*}
&\bigg| \int_{|v|\leqslant C} \nu F \int_0^\infty  P_0 \big[\tilde \varphi_\eps(t,x,v,z) - \varphi(t,x) - \varepsilon v \lceil v \rfloor^{-\beta} z \cdot \nabla_x\varphi(t,x) \big] \ud z \ud v \bigg| \\
&\leqslant \nu_2 \int_{|v| \leqslant C} \lceil v \rfloor^{\beta} F \int_0^\infty\varepsilon ^2 e^{-\nu_1 z} |v|^{2} \lceil v \rfloor^{-2\beta} z^2\| \na_x^2\varphi\|_{L^{\infty}} \ud z\ud v \\
&\leqslant C\varepsilon^2\| \na_x^2\varphi\|_{L^{\infty}} .
\end{align*}
Hence, the estimate and limit of  Lemma \ref{petites vitesses BL} hold for $\alpha \geqslant 1$. For the case $\alpha < 1$, the estimate follows directly from expression \eqref{petites vitesses alpha<1} together with inequality \eqref{estimation phi_eps -phi} for $|v| \leqslant C$, and the limit is a direct consequence of the inequality since in this case $1-\gamma = \frac{1-\alpha}{1-\beta} > 0$. \epl \\

\ni \textbf{Step 2: Convergence to the elliptic operator} \\
We still have to deal with the following limit
\begin{equation}\label{lim L-eps}
\underset{\eps \to 0}{\lim} \ L^\eps_\varphi := \underset{\eps \to 0}{\lim} \ \eps^{-\gamma} \int_{|v| \geqslant C}  \nu F\left(\chi^{\varepsilon}-\varphi - \frac{\eps}{\nu} \ j^\eps_F \cdot \na_x\varphi \right) \ud v .
\end{equation}
For $|v| \geqslant C$,  $\lceil v \rfloor^{-\beta} = |v|^{-\beta}$.  By performing the change of variable $ w = \varepsilon z|v|^{-\beta}v$, we obtain
$$
|w| = \varepsilon z|v|^{1-\beta} ,  \ \ |v| = \big(\frac{|w|}{\varepsilon z}\big)^\frac{1}{1-\beta} ,  \  v = \frac{w}{(\varepsilon z)^{\frac{1}{1-\beta}}|w|^{-\frac{\beta}{1-\beta}}}  \ \mbox{ and } \ \ud v = \frac{1}{1-\beta}\frac{\ud w}{(\varepsilon z)^{\frac{d}{1-\beta}}|w|^{-\frac{\beta d}{1-\beta}}} . $$
Thus, the integral of limit \eqref{lim L-eps} becomes \\
\textbf{Case 1: $\alpha <1$.}
\begin{align*}
L^\eps_\varphi &: = \eps^{-\gamma} \int_{|v| \geqslant C} \nu(x,v) F(v) \int_0^\infty P_{\eps}(x,v,z) \big[\varphi(t,x+\eps |v|^{-\beta}vz) - \varphi(t,x)\big] \ud z \ud v \\
&= \frac{\kappa}{1-\beta} \ \int_0^\infty \int_{|w| \geqslant C \varepsilon z} z^{\gamma} e^{-z \int_0^1 \nu_0(x+sw) \ud s} \ \nu_0(x)\nu_0(x+w) \ \frac{\varphi(t,x+w) - \varphi(t,x)}{|w|^{\gamma +d}} \ \ud w\ud z .
\end{align*}
\textbf{Case 2: $\alpha = 1$.}
\begin{align*}
L^\eps_\varphi &: = \eps^{-\gamma} \int_{|v| \geqslant C} \nu(x,v) F(v) \int_0^\infty P_{\varepsilon}(x,v,z) [\tilde \varphi_\eps(t,x,v,z) -\varphi(t,x)] \ud z \ud v \\
& \qquad - \eps^{-\gamma} \int_{C\leqslant |v| \leqslant \eps^{-\frac{1}{1-\beta}}} \nu(x,v) F(v) \int_0^\infty P_0(x,v,z) \ \eps |v|^{-\beta}vz \cdot \na_x \varphi(t,x) \ud z \ud v \\
&= \frac{\kappa}{1-\beta} \ \int_0^\infty \int_{|w| \geqslant C \varepsilon z} z^{\gamma} e^{-z \int_0^1 \nu_0(x+sw) \ud s} \ \nu_0(x)\nu_0(x+w) \ \frac{\varphi(t,x+w) - \varphi(t,x)}{|w|^{\gamma +d}} \ \ud w\ud z \\
& \qquad -\frac{\kappa}{1-\beta} \ \int_0^\infty \int_{C \eps z \leqslant |w| \leqslant 1} z^{\gamma} e^{-z \nu_0(x) } \ \nu_0(x)\nu_0(x) \ \frac{w \cdot \na_x \varphi(t,x)}{|w|^{\gamma +d}} \ \ud w\ud z .
\end{align*}
\textbf{Case 3: $\alpha >1$.}
\begin{align*}
L^\eps_\varphi &: = \eps^{-\gamma} \int_{|v| \geqslant C} \nu(x,v) F(v) \int_0^\infty P_{\varepsilon}(x,v,z) [\tilde \varphi_\eps(t,x,v,z) -\varphi(t,x)] \ud z \ud v \\
& \qquad - \eps^{-\gamma} \int_{|v| \geqslant C} \nu(x,v) F(v) \int_0^\infty P_0(x,v,z) \ \eps |v|^{-\beta}vz \cdot \na_x \varphi(t,x) \ud z \ud v \\
&= \frac{\kappa}{1-\beta} \ \int_0^\infty \int_{|w| \geqslant C \varepsilon z} z^{\gamma} e^{-z \int_0^1 \nu_0(x+sw) \ud s} \ \nu_0(x)\nu_0(x+w) \ \frac{\varphi(t,x+w) - \varphi(t,x)}{|w|^{\gamma +d}} \ \ud w\ud z \\
& \qquad -\frac{\kappa}{1-\beta} \ \int_0^\infty \int_{|w| \geqslant C \varepsilon z} z^{\gamma} e^{-z \nu_0(x) } \ \nu_0(x)\nu_0(x) \ \frac{w \cdot \na_x \varphi(t,x)}{|w|^{\gamma +d}} \ \ud w\ud z .
\end{align*}
Hence, in all three cases,
\begin{equation}\label{lim L-eps = kappa L}
	\lim_{\eps \to 0} \ L^\eps_\varphi  = \frac{\kappa}{1-\beta} \ \mathrm{PV} \int_{\RR^{d}} \eta(x,y)\frac{\varphi(t,y) - \varphi(t,x)}{|x-y|^{\gamma+d}} \ \ud y  =: -\kappa \ \mathcal{L}^*(\varphi) 
\end{equation}
with 
$$ \eta(x,y) = \nu_0(x)\nu_0(y) \int_0^\infty z^{\gamma}  e^{ - z\int_0^1 \nu_0(sx+(1-s)y)\ud s} \ \ud z .$$
\begin{remark}
Note that the integral corresponding to the term containing $w|w|^{-\gamma - d} \cdot \na_x \varphi$ is zero due to symmetry in $w$. We keep it here to justify that the integrals are well-defined near the origin $0$. Thus, we have a unified formula for $L^\eps_\varphi$ for all $\alpha >0$.
\end{remark}
In order to justify limit \eqref{lim L-eps = kappa L} rigorously and show that it holds uniformly with respect to $x$ and $t$, we introduce the function
$$ \bar\eta(x,y,z) = \frac{\kappa}{1-\beta} \ \nu_0(x)\nu_0(y) z^{\gamma}  e^{- z\int_0^1 \nu_0(sx+(1-s)y) \ud s} $$ 
and we split the integral $L^\eps_\varphi$ as follows: 
\begin{align*}
L^\eps_\varphi = &\int_0^\infty \int_{|w| \geqslant C \varepsilon z} \bar\eta(x,x+w,z) \ \frac{\varphi(t,x+w) - \varphi(t,x)}{|w|^{\gamma +d}}\ud w\ud z \\
= &\int_0^{\frac{1}{C\varepsilon}} \int_{|w| \geqslant 1} \bar\eta(x,x+w,z) \ \frac{\varphi(t,x+w) - \varphi(t,x)}{|w|^{\gamma +d}}\ud w\ud z \\
&+ \int_0^{\frac{1}{C\varepsilon}} \int_{C\varepsilon z \leqslant |w| \leqslant 1} \big[ \bar\eta(x,x+w,z) -\bar\eta(x,x,z)\big] \ \frac{\varphi(t,x+w) - \varphi(t,x)}{|w|^{\gamma +d}} \ \ud w \ud z \\
&+ \int_0^{\frac{1}{C\varepsilon}} \int_{C\varepsilon z \leqslant |w| \leqslant 1} \bar\eta(x,x,z) \ \frac{\varphi(t,x+w) - \varphi(t,x) - w \cdot\nabla_x\varphi}{|w|^{\gamma +d}} \ \ud w \ud z \\
&+ \int_{\frac{1}{C \varepsilon}}^\infty \int_{|w| \geqslant C \varepsilon z} \bar\eta(x,x+w,z) \ \frac{\varphi(t,x+w) - \varphi(t,x)}{|w|^{\gamma +d}} \ \ud w\ud z .
\end{align*}
The function $\bar \eta$ plays the role of $P_\eps$ for small velocities and was introduced in order to make the previous quantities integrable.  Note that all the integrals above are defined in the classical sense without need for principal value and the fact that we have assumed that $\nu(x,v) = \nu_0(x)|v|^{-\beta}$ and $F(v)=C |v|^{-\alpha-d}$ for $|v| \geqslant C$, made passing to the limit much easier ($\lim_{\varepsilon\rightarrow 0} \ L^\eps_\varphi$ is exactly the operator defined by the principal value).    
 This completes the proof of Theorem \ref{main BL}.
 
 \section*{Some comments}
We conclude this paper with some comments on the method used.
 \begin{itemize}
 	\item Although this method is simple and robust for handling the case of spatial dependence, it still relies on the structure of the operator. Indeed, despite the absence of a spectral gap caused by the fact that the collision frequency $\nu$ is unbounded and may degenerate for large velocities, we were able to take the limit in the weak formulation, specifically in the term corresponding to the projection of the solution onto the kernel. This was achieved thanks to the inequality $|Q^+(g)| \lesssim \nu F \|g\|_{L^2_{\nu F^{-1}}}$, which follows from the integral structure of $Q$.  
\item This method is very difficult to adapt in the case of the Fokker-Planck operator with a heavy-tailed equilibrium, due to the absence of the spectral gap and the estimate of the type mentioned in the previous point. In addition to this, the decomposition of the operator into $Q^+$ and $Q^-$ is not possible in the sense that the terms cannot be treated separately as in Proposition \ref{partie negligeable}. At some level, all terms of the operator, including the advection term, are of the same order.	
 \end{itemize}

{\footnotesize
 \bibliography{bib-BL}
}

\end{document}